\documentclass[final,3p,numbers,sort&compress]{elsarticle}

\usepackage{amsmath}
\usepackage{amssymb}

\usepackage{color}
\usepackage{array}
\usepackage{booktabs}
\usepackage{enumitem}
\usepackage[misc,geometry]{ifsym}

\newcommand{\bs}{\boldsymbol}

\DeclareMathAlphabet\mathbfcal{OMS}{cmsy}{b}{n}

\usepackage{graphicx}
\usepackage{mathrsfs}
\usepackage{color}
\usepackage[bookmarks,bookmarksnumbered]{hyperref}
\hypersetup{colorlinks = true,
linkcolor = blue,
anchorcolor =red,
citecolor = blue,
filecolor = red,
urlcolor = red,
pdfauthor=author}

\usepackage{tikz}
\usetikzlibrary{calc,positioning}
\newtheorem{remark}{Remark}

\journal{Physica D: Nonlinear Phenomena}
\begin{document}

\begin{frontmatter}
\title{Model reduction of nonlinear time-delay systems via ODE approximation and spectral submanifolds}
\author[label1]{Yuan Tang}
\author[label1]{Mingwu Li\corref{cor1}}
\ead{limw@sustech.edu.cn}
\cortext[cor1]{Corresponding author}
\address[label1]{Department of Mechanics and Aerospace Engineering, Southern University of Science and Technology, Shenzhen, China}

\begin{center}
  \textbf{(Accepted by \textit{Physica D})}
\end{center}

\begin{abstract}
Time-delay dynamical systems inherently embody infinite-dimensional dynamics, thereby amplifying their complexity. This aspect is especially notable in nonlinear dynamical systems, which frequently defy analytical solutions and necessitate approximations or numerical methods. These requirements present considerable challenges for the real-time simulation and analysis of their nonlinear dynamics. To address these challenges, we present a model reduction framework for nonlinear time-delay systems using spectral submanifolds (SSMs). We first approximate the time-delay systems as ordinary differential equations (ODEs) without delay and then compute the SSMs and their associated reduced-order models (ROMs) of the ODE approximations. These SSM-based ROMs successfully predict the nonlinear dynamical behaviors of the time-delay systems, including free and forced vibrations, and accurately identify critical features such as isolated branches in the forced response curves and bifurcations of periodic and quasi-periodic orbits. The efficiency and accuracy of the ROMs are demonstrated through examples of increasing complexity.
\end{abstract}

\begin{keyword}
Nonlinear time-delay systems \sep Invariant manifolds \sep Reduced-order models \sep Spectral submanifolds \sep Bifurcation
\end{keyword}
\end{frontmatter}

\section{Introduction}
Time-delay phenomena profoundly impact the system’s dynamic behavior, and even the simplest time-delay dynamical systems can exhibit extremely complex behavior patterns~\cite{gu2003stability,reddy2000experimental,beta2004controlling}. Indeed, a small delay can significantly impact the stability and control performance of the entire controlled system~\cite{hu2003dynamics}. It can also affect synchronous oscillations and stability switches, causing amplitude death of coupled neurons~\cite{reddy2000experimental}. The role of time delays is not limited to negative effects. They can be utilized to enhance a system's performance. Time-delayed feedback control has appeared as an effective method to suppress chaotic motion~\cite{ding2023control,pyragas1992continuous,pyragas1993experimental,ikeda1987high} and stabilize unstable states~\cite{hu2003dynamics,campbell2007time}.

Time-delay systems are characterized by delay differential equations (DDEs). The evolution of a DDE's state depends not only on the current state but also on the historical state at some previous moment or over a certain duration, which highlights the fact that systems with time delays are infinite-dimensional~\cite{gu2003stability}. Various numerical methods are developed and used to discretize and approximate these time-delay systems, resulting in high-dimensional nonlinear systems that are associated with high computational costs and substantial memory demands~\cite{gyori1991uniform, Westdal1970,chavez2020numerical}.

To address these challenges, various model reduction techniques have been proposed to streamline the analysis and control of high-dimensional nonlinear systems caused by time delays. For linear time-delay systems, projection of the original system onto some carefully selected subspace results in low-dimensional, reduced-order models. Following this line, Krylov subspace methods and balanced truncation techniques are commonly used for performing model reduction of DDEs~\cite{michiels2011krylov}. Balanced truncation achieved via solving Lyapunov equations can construct a reduced-order model preserving systems' stability and performance~\cite{Jarlebring2013Model}. As an extension of existing balanced truncation reduction techniques, a delay-dependent approach based on energy functionals that characterize the observability and controllability properties of the time delay system was recently proposed~\cite{lordejani2020model}. Expansion via orthogonal Laguerre basis functions has also been shown as an effective reduction scheme~\cite{samuel2014model,wang2016laguerre}.

Model reduction for nonlinear time-delay systems is much more challenging. A nonlinear time-delay system can be approximated as piecewise linear time-delay systems, which enables the application of linear reduction techniques~\cite{wang2013stability}. However, this approximation is problem-dependent and not suitable for highly nonlinear systems. For nonlinear DDEs with periodic coefficients, the Chebyshev spectral collocation method was used to achieve dimensional reduction~\cite{deshmukh2008dimensional}. For DDEs with a locally exponentially stable equilibrium, moments of linear systems can be extended to nonlinear systems via center manifold theory~\cite{scarciotti2014model,scarciotti2015model}. Such moments characterize the steady-state response of the output of an interconnected system composed of the DDEs and a signal generator. Then, reduced-order models that satisfy moment matching can be constructed~\cite{scarciotti2014model,scarciotti2015model}.

Model reduction methods based on invariant manifolds have received considerable attention in the past two decades~\cite{pesheck2000reduced,pesheck2002new,jiang2004nonlinear,touze2006nonlinear,haller2016nonlinear,touze2021model}. An exact reduction can be obtained using the invariant manifolds, which predicts precisely nonlinear and long-term dynamics. This is crucial for designing effective control strategies and understanding the fundamental physics of the DDEs. Reduction on center manifolds has been used to construct ROMs for nonlinear time-delay systems composed of rigid and flexible substructures~\cite{Wang2001Dimensional}. Recently, nonlinear normal modes~\cite{Shaw1999ModalAR, Kerschen2009Nonlinear} were used to build ROMs for delay-coupled limit cycle oscillators~\cite{govind2023nonlinear}. The constructions in~\cite{govind2023nonlinear} assumed a small delay and were truncated at cubic order.

Center manifolds are not structurally stable because they are not persistent under the addition of damping. In addition, they are not unique. Likewise, nonlinear normal modes are not unique even for linear dynamical systems~\cite{haller2016nonlinear}. In contrast, spectral submanifolds (SSM) defined for hyperbolic fixed points are structurally stable and also unique under appropriate non-resonance conditions~\cite{haller2016nonlinear,haller2023nonlinear}. Therefore, SSMs have emerged as a powerful framework for constructing low-dimensional ROMs for high-dimensional mechanical systems, both in equation-driven~\cite{jain2022compute,li2022nonlinear1,li2022nonlinear2,li2023model,li2024fast,li2023nonlinear} and data-driven settings~\cite{cenedese2022data,axaas2023fast,liu2024data}. In particular, SSM-based model reduction can be used to effectively predict backbone curves~\cite{szalai2017nonlinear,breunung2018explicit}, forced response curves~\cite{breunung2018explicit,jain2022compute,li2022nonlinear1,li2024fast}, quasi-periodic and chaotic motions~\cite{li2022nonlinear2,liu2024data}, and bifurcations~\cite{li2022nonlinear2,li2023nonlinear}.

An SSM is defined as the smoothest nonlinear continuation of a spectral subspace of a linear system under the addition of nonlinearity and possibly external forcing, which can be periodic, quasi-periodic, or even chaotic~\cite{haller2016nonlinear,haller2024nonlinear}. SSMs can be computed via parameterization methods~\cite{haro2016parameterization}. In earlier studies, equations of motion were assumed to be in modal coordinates, such that the coefficient matrix was diagonal. This assumption hinders the application of SSM reduction to high-dimensional finite element (FE) problems, where the transformation from physical coordinates to modal coordinates is infeasible for high-dimensional FE problems~\cite{jain2022compute}. To address this challenge, schemes that compute SSMs directly in physical coordinates were recently developed~\cite{jain2022compute,vizzaccaro2021direct}. In particular, the computational scheme in~\cite{jain2022compute} enables the computation of SSMs of arbitrary dimensions and order of accuracy for general dynamical systems. This computational scheme has been implemented in the open-source package SSMTool~\cite{jain2021ssmtool}.

In the above-reported studies, SSM reductions were mainly applied to systems whose equations of motion are in the form of ordinary differential equations (ODEs), with an exception to systems in the form of differential-algebraic equations (DAEs)~\cite{li2023model}. Here, we aim at constructing SSM-based ROMs for nonlinear time-delay dynamical systems described by DDEs. We note that SSM-based reduction has been recently applied to autonomous delay systems with nonlinearity up to cubic orders~\cite{SZAKSZ202467,szaksz2025spectral}. They formulated DDEs as operator differential equations and then parameterized the SSMs with delay-dependent expansion coefficients. The derivations in~\cite{SZAKSZ202467,szaksz2025spectral}, however, are truncated at cubic expansion, which might not be sufficient for highly nonlinear problems. In addition, the derivations in~\cite{SZAKSZ202467,szaksz2025spectral} are tailored to autonomous dynamics, and external forcing is not considered.

Here, we address the above limitations by focusing on an alternative methodology to leverage the existing developments of SSM reduction of ODEs fully.
Specifically, the infinite-dimensional time-delay systems are first approximated using high-dimensional ODEs. Then SSM-based ROMs are constructed for the ODEs using SSMTool that automates the computation of SSM-based ROMs up to arbitrary expansion orders. These ROMs will be utilized to predict the system’s nonlinear dynamics, including free vibration, self-excited oscillation, the isolas of forced response curves, bifurcations of periodic orbits, etc. We will consider several time-delay nonlinear systems ranging from simple to complex to demonstrate the effectiveness of the proposed reduction framework.

The rest of this paper is organized as follows. We present the setup of DDEs and problem statements in Sect.~\ref{sec: setup}, followed by a discussion on the transformation from infinite-dimensional DDEs to high-dimensional ODEs in Sect.~\ref{sec:dde-to-ode}. We then present SSM-based reduction for the transformed ODEs in Sect.~\ref{sec:ssm-reduction}, where we show that we can consistently reduce the high-dimensional ODEs to two-dimensional ROMs, independently of the dimensionality of the ODEs. We use a suite of examples with increasing complexity in Sect.~\ref{sec:example} to illustrate the effectiveness of the SSM-based reduction. We conclude this study in Sect.~\ref{sec:conclusion}.

\section{Setup and problem statement}
\label{sec: setup}
Consider a system of delay differential equations (DDEs) below.
\begin{equation}
\label{eq: DDEs_form}
\dot{\bs{x}}=\bs{f}(\bs{x},\bs{x}(t-\tau_d))+\epsilon\bs{g}(\Omega t),\quad \epsilon\ll 1,
\end{equation}
where \textbf{$\bs{x}(t) \in \mathbb{R}^n$} represents the current state of the system, $\tau_d > 0$ denotes a delay, and $\bs{f}(\bs{x}(t), \bs{x}(t-\tau_d))$ is the nonlinear function describing the system's autonomous dynamics. The term $\epsilon\mathbf{g}(\Omega t)$ represents an external harmonic excitation with forcing amplitude $\epsilon$ and frequency $\Omega$.

At $\epsilon=0$, we assume that the origin $\bs{x}=\bs{0}$ is a hyperbolic fixed point of the dynamical system, and we are interested in the solution to~\eqref{eq: DDEs_form} with the initial condition below
\begin{equation}
\label{eq: DDEs_initi}
\bs{x}(s) = \bs{x}_0(s),\quad s\in[-\tau_d,0],
\end{equation}
and also the limit cycle oscillations of~\eqref{eq: DDEs_form}. For such a limit cycle, it satisfies the periodic boundary conditions: $\bs{x}(t)=\bs{x}(t+T)$, where $T$ is an unknown period of the limit cycle.

When $\epsilon>0$, we seek periodic motions of the system under the external harmonic forcing. We are also interested in the bifurcations of the forced periodic orbits under variations in the forcing amplitude $\epsilon$ and frequency $\Omega$. In particular, when the periodic orbit undergoes a secondary Hopf bifurcation, a family of tori on which quasi-periodic orbits stay appears. We aim to extract these quasi-periodic orbits and further infer the bifurcations of these quasi-periodic orbits.

To solve these problems, we first transform the DDEs~\eqref{eq: DDEs_form} as systems of ordinary differential equations (ODEs) and then use the theory of spectral submanifold (SSM) to reduce these high-dimensional ODEs as low-dimensional reduced-order models (ROMs), which enables efficient and analytic predictions on the nonlinear dynamics of the original systems with delay. In the following two sections, we describe the transformation and the reduction in detail.

\begin{remark}
 We note that SSM theory was established for finite-dimensional ODEs~\cite{haller2016nonlinear}. An extension to infinite-dimensional analysis for a forced nonlinear beam model whose equations of motion are partial-differential equations was given in~\cite{kogelbauer2018rigorous}. These results are based on Cabré et al.~\cite{Cabre2002,Cabre2003,Cabre2005}, where relevant theorems are established in Banach spaces. It is out of the scope of this study to establish the theory of SSM for DDEs whose solution space is also a Banach space with a properly defined norm. Instead, we focus on the computational aspect since DDEs can be approximated via high-dimensional ODEs, as demonstrated in previous studies. We review the demonstration briefly in the next section.
\end{remark}

\section{Transforming DDEs as ODEs}
\label{sec:dde-to-ode}

One approach to studying the dynamics of DDEs is to approximate the original equations as high-dimensional systems of ordinary differential equations (ODEs). This approach was initially proposed by Repin, who used first-order approximations of the original solution to approximate DDEs with ODE systems, now known as the chain method~\cite{repin1965approximate}. Gyori and Turi demonstrated the uniform convergence of the chain method and discussed its applications to DDEs on infinite intervals ~\cite{gyori1991uniform}. Subsequently, the chain method has been applied to optimal time control ~\cite{Westdal1970,banks1979approximation}, linear differential equations with large time delays ~\cite{hess1972optimal}, and delay chemical reaction networks ~\cite{liptak2018approximation}. Recently, Chávez et al. proposed a new numerical method combining the chain method with a second-order approximation scheme using finite-order Taylor expansions. They showed that the new method can effectively approximate nonsmooth DDEs ~\cite{chavez2020numerical}.

Here, we use the numerical method proposed by Chávez et al.~\cite{chavez2020numerical} to perform the transformation. Assume that the solution to~\eqref{eq: DDEs_form} is smooth. Take $N\in\mathbb{N}$ sufficiently large and define grid points $\tau_{i}=i \tau_{d}/{N}$ for $i=0, \ldots, N$. Furthermore, we define $\bs{u}_i(t)=\bs{x}(t-\tau_i)$ for all $t\geq 0$ and $i=0,1,\ldots,N$. With Taylor expansion, we obtain
\begin{align}
    \bs{u}_{i-1}\left( t \right) & = \bs{x}\left( t-\tau_{i-1} \right)  = \bs{x}\left( t-\left( \tau_i-\frac{\tau_d}{N} \right) \right) = \bs{u}_i\left( t+\frac{\tau_d}{N} \right) =\sum_{k=0}^M{\frac{1}{k!}}\bs{u}_{i}^{\left( k \right)}\left( t \right) \left( \frac{\tau_d}{N} \right) ^k+\mathcal{O}\left( \left( \frac{\tau_d}{N} \right) ^{M+1} \right) \label{eq:DDEs_transform}
\end{align}
for $i=1,\ldots,N$. 
By selecting $M=2$ and introducing the auxiliary function $\bs{w}_i=\dot{\bs{u}}_i$, the aforementioned approximation transforms the original DDEs into a chain of equations, resulting in the following set of ODEs without delay
\begin{equation} 
\label{eq:dde-as-odes}
\left\{
\begin{aligned}
 \dot{\bs{u}}_0(t) &= \bs{f}( \bs{u}_0(t), \bs{u}_N(t))+\epsilon\bs{g}(\Omega t), & t &\geq 0, \\
 \dot{\bs{u}}_i(t) &= \bs{w}_i(t), & t &\geq 0, \quad i = 1, \ldots, N, \\
 \dot{\bs{w}}_i(t) &= \frac{2N^2}{\tau_d^2} \left( \bs{u}_{i-1}(t) - \bs{u}_i(t) - \frac{\tau_d}{N}\bs{w}_i(t) \right), & t &\geq 0, \quad i = 1, \ldots, N.
\end{aligned}
\right.
\end{equation}

\section{Reduction to spectral submanifolds}
\label{sec:ssm-reduction}
\subsection{Spectral submanifold and model reduction}
We separate the linear and nonlinear parts of the vector field $\bs{f}( \bs{u}_0, \bs{u}_N)$ in~\eqref{eq:dde-as-odes} as below
\begin{equation}
\label{eq:vector-eq}
 \bs{f}( \bs{u}_0, \bs{u}_N)=\boldsymbol{A}_{\boldsymbol{u}_0}\bs{u}_0+\boldsymbol{A}_{\boldsymbol{u}_N}\bs{u}_N+\boldsymbol{f}_\mathrm{nl}( \bs{u}_0, \bs{u}_N),
\end{equation}
where $\bs{A}_{\bs{u}_0},\bs{A}_{\bs{u}_N}\in\mathbb{R}^{n\times n}$ are coefficient matrices of linear terms, and $\boldsymbol{f}_\mathrm{nl}$ stands for a smooth nonlinear function. Let $\bs{z}=(\bs{u}_0,\bs{u}_1,\cdots,\bs{u}_N,\bs{w}_1,\cdots,\bs{w}_N)\in\mathbb{R}^{(2N+1)n}$, the dynamical system~\eqref{eq:dde-as-odes} can be rewritten in the standard form below:
\begin{equation}
\dot{\bs{z}}=\bs{A}\bs{z}+\bs{F}\left( \bs{z} \right) +\epsilon \bs{F}^{\text{ext}}\left( \Omega t \right),
\label{eq:non_eq1}
\end{equation}
where
\begin{equation}
 \bs{A}=\left(\begin{array}{cc} \bs{A}_{\boldsymbol{u}_0\boldsymbol{u}_0}&\bs{A}_{\boldsymbol{u}_0\boldsymbol{w}}\\
 \bs{A}_{\boldsymbol{u}_i\boldsymbol{u}_i}& \bs{A}_{\boldsymbol{u}_i\boldsymbol{w}}\\
 \bs{{A}_{wu}}&\bs{{A}_{ww}}\end{array}\right), \quad 
 \bs{F(z)} = \left(\begin{array}{cc}
 \bs{f}_\mathrm{nl}( \bs{u}_0, \bs{u}_N) \\
 \bs{0} \\
 \bs{0}\end{array}\right),\quad 
 \bs{F}^\mathrm{ext} = \left(\begin{array}{cc}
 \epsilon\bs{g}(\Omega t) \\
 \textbf{0} \\
 \textbf{0}\end{array}\right),
\end{equation}
with
\begin{equation}
\begin{aligned}
\bs{A}_{\bs{u}_0\bs{u}_0} &= (\bs{A}_{\bs{u}_0}, \bs{0}, \cdots,\bs{0}, \bs{A}_{\bs{u}_N})_{n\times (N+1)n},\quad 
\bs{A}_{\bs{u}_0\bs{w}} = \bs{0}_{n\times Nn},\\
\bs{A}_{\bs{u}_i\bs{u}_i} &= \bs{0}_{Nn\times (N+1)n},\quad 
\bs{A}_{\bs{u}_i\bs{w}} = \bs{I}_{Nn},\\
\bs{A}_{\bs{w}\bs{u}} &= \frac{2N^2}{\tau_d^2}\begin{bmatrix}\bs{I}_n &-\bs{I}_n \\ & \ddots & \ddots \\ & & \bs{I}_n& -\bs{I}_n\end{bmatrix}_{Nn\times (N+1)n}, \quad
\bs{A}_{\bs{w}\bs{w}} = -\frac{2N}{\tau_d}\begin{bmatrix}\bs{I}_n \\ & \ddots \\ & & \bs{I}_n
\end{bmatrix}_{Nn\times Nn},
\end{aligned}
\end{equation}
and $\bs{I}_{Nn}$ and $\bs{I}_n$ being identity matrices of dimension $Nn$ and $n$, respectively.

Let $\text{Spect}(\bs{A}) = \{\lambda_1, \ldots, \lambda_{2(N+1)n}\}$ be the spectrum of the linear part of the dynamical system. We have
\begin{equation}
\bs{Av}_i=\lambda _i\bs{v}_i,\quad 1\le i\le 2(N+1)n,
\end{equation}
where $\bs{v}_i$ is the eigenvector corresponding to the eigenvalue ${\lambda}_i$. In this paper, we aim at constructing two-dimensional ROMs to achieve a significant dimension reduction to the system~\eqref{eq:non_eq1}. This is indeed feasible for systems without internal resonances~\cite{haller2016nonlinear}. Internal resonances arise when specific integer relationships exist between the eigenvalues of a system, i.e., $n_1\lambda_1^\mathcal{E}+n_2\lambda_2^\mathcal{E} = \lambda_j$, where $\lambda_{1,2}^\mathcal{E}$ are the eigenvalues associated with the master spectral subspace $\mathcal{E}$ for reduction, $\lambda_j$ is any other eigenvalue, $n_1$ and $n_2$ are some non-negative integers. These resonances lead to energy transfer among modes and can necessitate higher-dimensional reductions~\cite{li2022nonlinear1}. In the absence of such resonances, as considered in this study, a two-dimensional SSM provides an accurate and computationally efficient reduced-order model. Here, we take $\mathcal{E}=\text{span}\left( v^{\mathcal{E}},\bar{v}^{\mathcal{E}} \right) $ as a master subspace for model reduction, where the two eigenvectors that span $\mathcal{E}$ correspond to a pair of complex conjugate eigenvalues $ \left( \lambda ^{\mathcal{E}},\bar{\lambda}^{\mathcal{E}} \right) $ that are in $\text{Spect}(\bs{A})$ . Specifically, we select the master subspace $\mathcal{E}$ such that the pair of complex conjugate eigenvalues$\left(\lambda ^{\mathcal{E}},\bar{\lambda}^{\mathcal{E}}\right)$ has the largest real parts. It follows that $\mathcal{E}$ consists of the slowest decaying modes if $\mathrm{Re}\left(\lambda ^{\varepsilon} \right) <0 $ or the fastest growing modes if $\mathrm{Re}\left(\lambda ^{\varepsilon}\right)>0$. However, if the dominant eigenvalue is real, the spectral subspace becomes one-dimensional, in which case the associated reduced-order model captures the dynamics along this one-dimensional manifold.

Under the addition of the nonlinear term $\bs{F}(\bs{z})$, the spectral subspace $\mathcal{E}$ is perturbed into some invariant manifolds of the system $\dot{\bs{z}}=\bs{A}\bs{z}+\bs{F}(\bs{z})$. It turns out that there are infinitely many such invariant manifolds that are tangent to $\mathcal{E}$ at the origin $\bs{z}=\bs{0}$. Among these invariant manifolds, however, there exists a unique, smoothest invariant manifold, which is defined as the SSM $\mathcal{W}(\mathcal{E})$ associated with the master subspace $\mathcal{E}$~\cite{haller2016nonlinear}. We note that the SSM is slow, attracting manifold if $\text{Re}(\lambda_{\mathcal{E}}) < 0$. Such a slow SSM lays the foundation for performing exact model reductions. Likewise, if $\text{Re}(\lambda_{\mathcal{E}}) >0$, the corresponding unstable SSM with underlying fast-growing modes provides a natural choice for model reduction.

The SSM $\mathcal{W}(\mathcal{E})$ can be parameterized via a map $\bs{z}=\bs{W}(\bs{p})$, where $\bs{p} = (p,\bar{p}) \in \mathbb{C}^2$ denotes a vector of parameterization coordinates. Meanwhile, the reduced dynamics on the SSM can be expressed as $\dot{\bs{p}} = \bs{R}(\bs{p})$. We can solve for the map $\bs{W}(\bs{p})$ and the vector field $\bs{R}(\bs{p})$ from the invariance equations associated with the SSM. We note that the solution to them can be obtained from SSMTool~\cite{jain2021ssmtool}, an open-source package that automates the computation of the SSM and the reduced dynamics up to any expansion orders.
Under the assumption 
the normal-form-style parameterization of the SSM gives the reduced dynamics in polar coordinates $\bs{p} = (\rho e^{\mathrm{i}\theta}, \rho e^{-\mathrm{i}\theta})$ below~\cite{breunung2018explicit}
\begin{equation}
\label{eq:red-auto}
\left\{
\begin{aligned}
\dot{\rho} &= a(\rho) = \rho \operatorname{Re}(\lambda^{\varepsilon}) + \sum_{j\geq 1}\operatorname{Re}(\gamma_j) \rho^{2j+1}, \\
\dot{\theta} &= b(\rho) = \operatorname{Im}(\lambda^{\varepsilon}) + \sum_{j\geq 1}\operatorname{Im}(\gamma_j) \rho^{2j},
\end{aligned}
\right.
\end{equation}
where $\{\gamma_j\}$ for $j \geq 1$ are some coefficients obtained from the parameterization. The vector field of the reduced dynamics~\eqref{eq:red-auto} is in the normal form, which is the simplest form of the reduced dynamics, ensuring that the nonlinear behavior of the original system can be accurately captured. 

Under further addition of the external harmonic excitation $\epsilon\bs{F}^\mathrm{ext}(\Omega t)$, the fixed point $\bs{z}=\bs{0}$ is perturbed as a periodic orbit $\gamma_\epsilon$, and accordingly, the autonomous SSM $\mathcal{W}(\mathcal{E})$ is perturbed into a time-periodic SSM $\mathcal{W}(\mathcal{E}, \Omega t)$ attached to the periodic orbit $\gamma_\epsilon$. In this case, the SSM parameterization map is updated as~\cite{breunung2018explicit,jain2022compute}
\begin{equation}
\bs{z} = \bs{W}(\bs{p})+\epsilon \bs{X}(\bs{p},\Omega t)+\mathcal{O}(\epsilon^2),
\end{equation}
and the reduced dynamic on the SSM is non-autonomous and can be expressed as~\cite{breunung2018explicit,jain2022compute}
\begin{equation}
\dot{\bs{p}}=\bs{R}\left(\bs{p}\right) +\epsilon\bs{S}\left(\bs{p},\Omega t \right) +\mathcal{O}\left( \epsilon ^2 \right).
\end{equation}
With a coordinate transformation $\theta = \theta + \Omega t$ and a leading-order approximation to the non-autonomous part $\bs{S}(\mathbf{p}, \Omega t)$, the reduced dynamics~\eqref{eq:red-auto} is updated as~\cite{jain2022compute}
\begin{equation}
\label{eq:red-non-auto}
\left\{
\begin{aligned}
\dot{\rho }&= a(\rho) + \epsilon \Omega^2(\operatorname{Re}(f) \cos \theta + \operatorname{Im}(f) \sin \theta),\\
\dot{\theta} &= b(\rho) - \Omega + \frac{\epsilon \Omega^2}{\rho}(\operatorname{Im}(f) \cos \theta - \operatorname{Re}(f) \sin \theta),
\end{aligned}
\right.
\end{equation}
where $f$ stands for a modal force.It is computed as the projection of the external forcing vector $\boldsymbol{F}^\mathrm{ext}$ in~\eqref{eq:non_eq1} onto the eigenvector $\boldsymbol{v}^{\mathcal{E}}$ associated with the master subspace. Interested readers can refer to~\cite{li2022nonlinear1} for detailed expressions of $f$. We note this modal force is distinct from the force vector \( \bs{f}(\bs{u}_0, \bs{u}_N) \) in~\eqref{eq:vector-eq}, which describes the system's autonomous dynamics, including both linear and nonlinear terms.  Here, we use SSMTool to obtain the ROM shown in~\eqref{eq:red-non-auto}.

\subsection{Predictions via SSM-based ROMs}
\label{sec:ssm-prediction}
The two-dimensional SSM-based ROM~\eqref{eq:red-auto} enables analytic predictions on the autonomous dynamics of the original time-delay system~\eqref{eq: DDEs_form}. Firstly, we can use the ROM to predict the system's transient responses. Further, $\Omega = b(\rho)$ gives a damped backbone curve of the ROM, which is further mapped as the backbone curve of the original system~\cite{breunung2018explicit}. Moreover, we can infer the existence of limit cycles on the SSM from~\eqref{eq:red-auto}. Specifically, suppose the polynomial $a(\rho)$ persistently admits a nontrivial, real root $\rho^\ast$ as the expansion order is increased, this root corresponds to a limit cycle with amplitude $\rho^\ast$ and period $2\pi/ b(\rho^\ast)$~\cite{ponsioen2019analytic,li2023nonlinear}. This limit cycle can be further mapped to the limit cycle of the time-delay system~\eqref{eq: DDEs_form}.

We can also make efficient and analytic predictions on the forced responses of the time-delay system~\eqref{eq: DDEs_form} using the ROM~\eqref{eq:red-non-auto}. Because of the transformation $\theta = \theta + \Omega t$, a fixed point of the vector field~\eqref{eq:red-non-auto} corresponds to a periodic orbit on the SSM of the original system. Moreover, the stability type of the fixed point is carried over to the stability type of the forced periodic orbit~\cite{li2022nonlinear1}. Therefore, one can make analytic predictions of the forced response curve of the system using the ROM~\eqref{eq:red-non-auto}, which enables effective extraction of isolated branches of the forced response curve~\cite{ponsioen2019analytic}. Further, one can detect the bifurcations of the forced periodic orbits of the original time-delay system~\eqref{eq: DDEs_form} as the bifurcations of the fixed points of the two-dimensional ROM, which significantly simplifies the demanding task of bifurcation detection of the infinite-dimensional time-delay system. Finally, we note that a limit cycle of the ROM~\eqref{eq:red-non-auto} corresponds to a torus on which a quasi-periodic orbit of the original time-delay system stays~\cite{li2022nonlinear2}. Thus, we can also make efficient predictions about the quasi-periodic orbits of the time-delay system~\eqref{eq: DDEs_form}.
 
In short, the reduction via SSMs enables us to address the challenges of high dimensionality of nonlinear time-delay systems. We will demonstrate it through several numerical examples in the following section.

\begin{remark}
\label{rmk:ssm-reduction}
    The above two-dimensional SSM-based reduction is effective provided that the system admits smooth nonlinearity, it has no internal resonance, and the response is within the domain of convergence of Taylor expansion of the SSM parameterization. We note that one can determine the domain of convergence via checking the convergence of backbone curves of $\Omega=b(\rho)$ in~\eqref{eq:red-auto}~\cite{li2023model}, or the roots of $a(\rho)$ in~\eqref{eq:red-auto}~\cite{ponsioen2019analytic}. We also note that the forcing amplitude $\epsilon$ should not be very large and the external forcing should not be in the form of parametric excitation because only leading-order approximation is adopted for the non-autonomous part of SSM in~\eqref{eq:red-non-auto}. In practice, we check the convergence of forced response curves to ensure its prediction accuracy~\cite{li2022nonlinear1}. Interested readers can refer to~\cite{thurnher2023nonautonomous} for two-dimensional SSM-based reduction for systems with parametric excitations.
\end{remark}

\section{Examples}
\label{sec:example}
In this section, we employ the high-dimensional ODE approximation approach to analyze infinite-dimensional time-delay systems, which has been detailed in Section~\ref{sec:dde-to-ode}. Subsequently, with the aid of the open-source package SSMTool~\cite{jain2021ssmtool}, we construct reduced-order models (ROMs) and use these ROMs to predict the system’s nonlinear dynamic behavior, encompassing phenomena such as free vibration, forced vibration response curves with isolas, and bifurcations of periodic orbits.

\subsection{A delayed Duffing oscillator}
\label{sec:ex-duffing}
We analyze the free and forced vibrations of a delayed Duffing oscillator:
\begin{equation}
\label{eq:duffing}
\ddot{x}=-\delta \dot{x}\left( t-\tau_d \right) -\alpha x-\beta x^3+\epsilon \cos \left( \Omega t \right).
\end{equation}
Let $x_1=x$ and $x_2=\dot{x}$, the equation above can be rewritten in the first-order form~\eqref{eq: DDEs_form} below
\begin{equation}
\label{eq:dde-duffing}
 \dot{x}_1=x_2,\quad \dot{x}_2=-\delta x_2(t-\tau_d)-\alpha x_1-\beta x_1^3+\epsilon\cos\Omega t.
\end{equation}
In the following computations, parameters are chosen as $\epsilon=0.0,\Omega=1.2,\alpha=2.0,\delta=0.2,\beta=-4.0,\tau_d=1.1$ unless otherwise stated.

We first follow Section~\ref{sec:dde-to-ode} to transform the original DDEs~\eqref{eq:dde-duffing} into a system of ODEs~\eqref{eq:non_eq1}. Recall that $N$ denotes the number of discrete points, we have $\bs{z}\in\mathbb{R}^{2(2N+1)}$ in this example. As detailed in~\ref{sec:app-ex1}, $N=100$ is sufficient to yield converged solutions. So we take $N=100$ in this example. Moreover, the analysis in~\ref{sec:app-ex1} shows that the system undergoes a Hopf bifurcation at $\tau_d=\tau^\ast=1.035$, as seen in the right panel of Fig.~\ref{fig:model1_r_distributed}.

\subsubsection{Reduction for dynamics before bifurcation}
We set $\tau_d=1.0<\tau^*$ and investigate the dynamics of the Duffing oscillator before the Hopf bifurcation. The first five eigenvalues in ascending order of their real parts for system~\eqref{eq:dde-duffing} are computed and listed below.\begin{equation}\label{eq:lamd-duffing-before-hopf}\lambda_{1,2}=-0.0057\pm 1.5182\mathrm{i},\quad\lambda_3=-2.8846+0.00\mathrm{i}, \quad\lambda_{4,5}=-3.6912\pm 7.3630\mathrm{i}.\end{equation}Since there is no internal resonance, we calculate the SSM tangent to the spectral subspace corresponding to the first pair of eigenvalues above, along with the reduced dynamics on the SSM. Thus, our SSM-based ROM is two-dimensional, independent of the dimension of the full-phase space.
First, we analyze the free vibration of system~\eqref{eq:duffing} with $\epsilon= 0$. Using SSMTool, we compute the two-dimensional slow SSM associated with $\lambda_{1,2}$, and the associated SSM-based ROM up to $\mathcal{O}(9)$ is as follows (cf.~\eqref{eq:red-auto}):
\begin{align}
\label{eq:ROM1_stable}
\dot{\rho} &= -1.45 \times 10^{-10}\rho^9 -1.576 \times 10^{-8}\rho^7 -2.099 \times 10^{-6}\rho^5 -0.0004142\rho^3 -0.005656\rho,\nonumber\\
\dot{\theta} &= -1.103\times 10^{-9}\rho^8 -1.261\times 10^{-7}\rho^6 -1.73\times 10^{-5}\rho^4 -0.003819\rho^2 +1.518.
\end{align}

\begin{figure}[!ht]
\centering
\includegraphics[width=.45\textwidth]{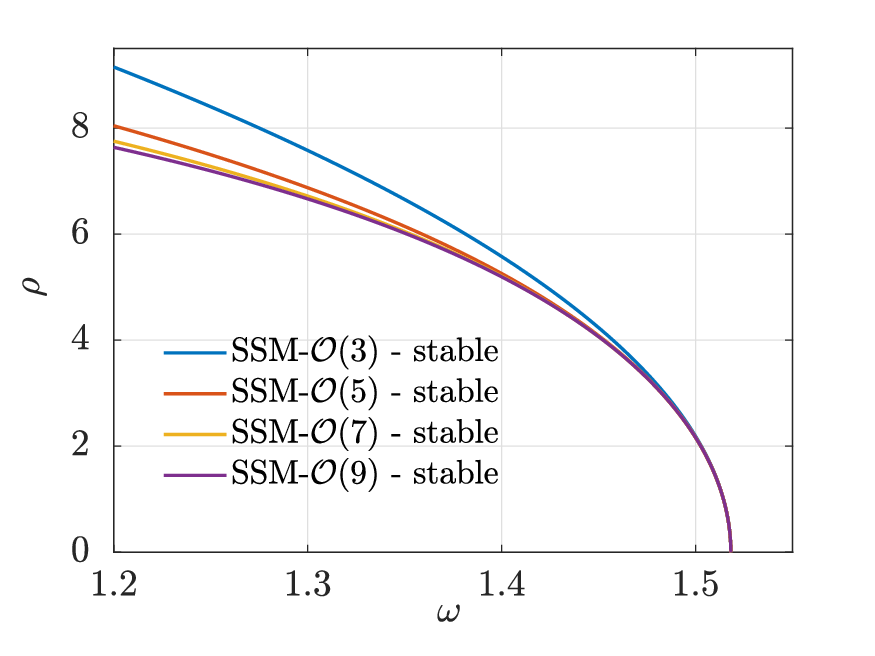}
\includegraphics[width=.45\textwidth]{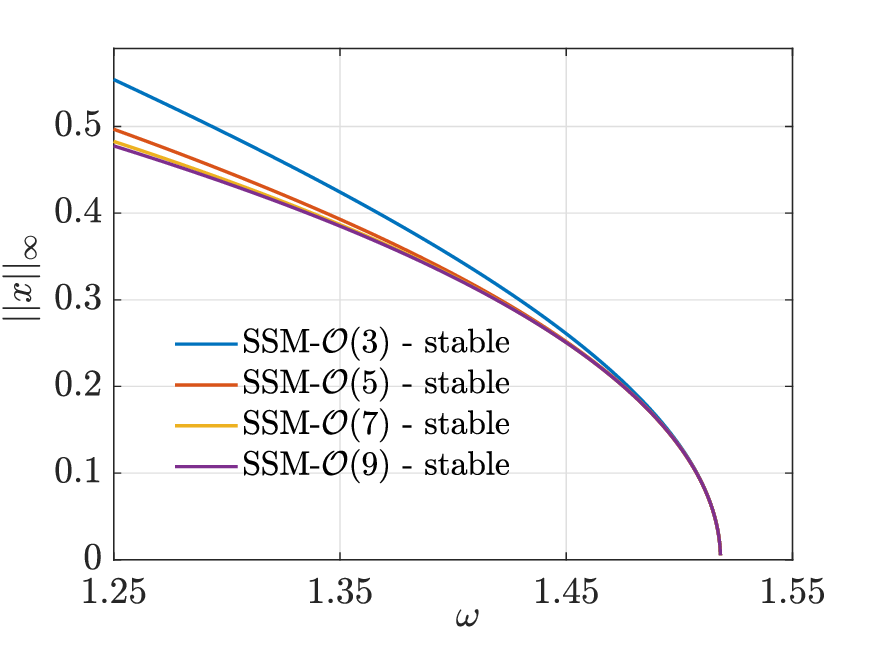}
\caption{\small Backbone curves in polar (reduced) and physical coordinates under increasing orders of expansion for the slowest two-dimensional SSM of the system~\eqref{eq:duffing} with $\tau_d=1.0$.}
\label{fig:model1_backbone_curve}
\end{figure}

The second sub-equation in the ROM determines the instantaneous oscillation frequency $\dot{\theta} = \omega$ as a function of the polar amplitude $\rho$, thus defining the damped backbone curve in polar coordinates~\cite{breunung2018explicit}. Fig.~\ref{fig:model1_backbone_curve} presents the backbone curve approximated at various expansion orders in both reduced and physical coordinates. It can be seen that the backbone curve tends to converge at higher amplitudes for SSM expansions at higher orders. For example, as seen in the left panel of Fig.~\ref{fig:model1_backbone_curve}, the backbone curve converges well with approximations of $\mathcal{O}(5)$ for $\rho \leq 3.75$, $\mathcal{O}(7)$ for $\rho \leq 5.85$, and $\mathcal{O}(9)$ for $\rho \leq 7.05$.

\begin{figure}[!ht]
\centering
\includegraphics[width=.45\textwidth]{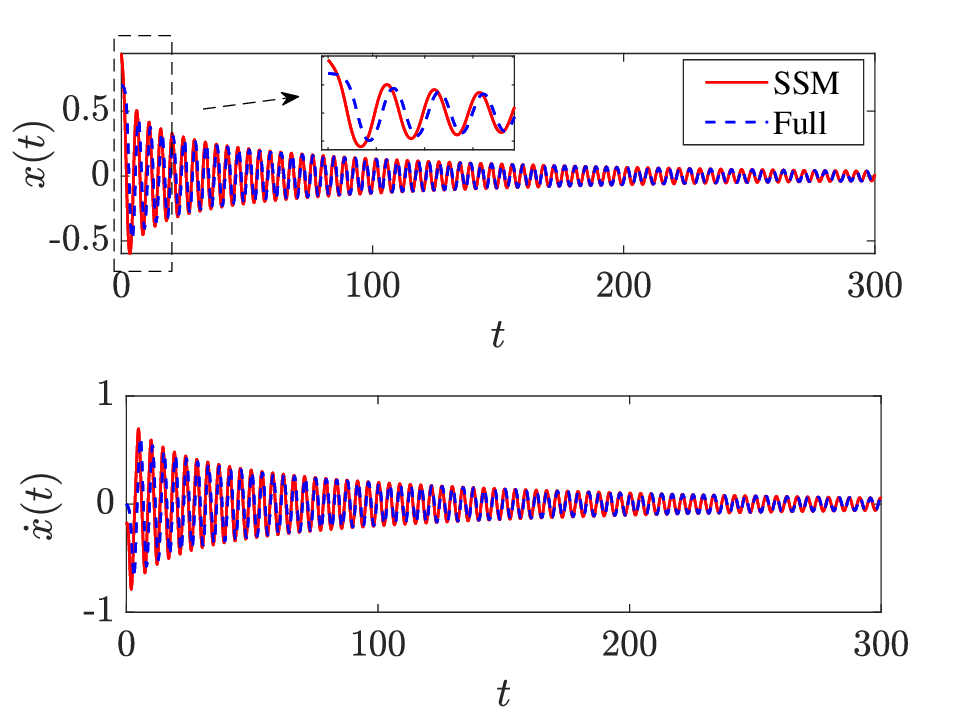} 
\includegraphics[width=.45\textwidth]{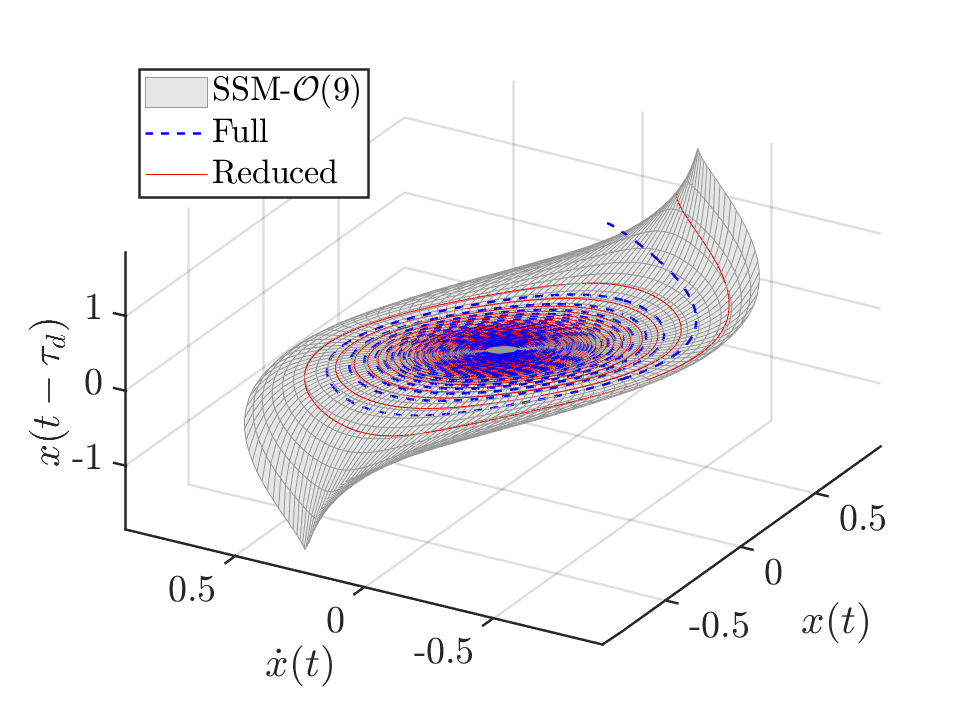}
\caption{\small (Left panel) Time history plots of \(x(t)\) and \(\dot{x}(t)\) for the Duffing oscillator~\eqref{eq:dde-duffing}. We note that the initial conditions~\eqref{eq: model1_dde_initial} are off the SSM, which explains the observed discrepancies at the initial stage. Such discrepancies disappear as time increases because the SSM is attracting. (Right panel) Visualization of the SSM projected onto the coordinates \((\dot{x}(t), x(t), \dot{x}(t-\tau_d))\) (gray plane) along with the trajectories of the full system (blue dashed line) and the prediction via SSM-based ROM~\eqref{eq:ROM1_stable} (red solid line).}
\label{fig:mode1_timehistory_pre}
\end{figure}

We can also use the SSM-based ROM~\eqref{eq:ROM1_stable} to make predictions on the free vibration of the system under given initial conditions. With initial conditions below
\begin{equation}
\label{eq: model1_dde_initial}
 x(s) = 0.7,\quad \dot{x}(s)= 0, \quad s\in[-\tau_d,0],
\end{equation}

 we perform forward simulation for the full system~\eqref{eq:non_eq1} and the SSM-based ROM~\eqref{eq:ROM1_stable}. The initial conditions for the SSM-based ROM are obtained following the transformation detailed in~\ref{sec:transform-ini}. The obtained results are plotted in Fig.~\ref{fig:mode1_timehistory_pre}, from which we see that the trajectory of the full system approaches the SSM very quickly and then synchronizes with the trajectory predicted from the SSM-based reduction.

In summary, the SSM-based ROM~\eqref{eq:ROM1_stable} can accurately capture the pre-bifurcation dynamics of the system~\eqref{eq:dde-duffing}. Moreover, making predictions via the SSM-based ROM~\eqref{eq:ROM1_stable} can greatly reduce the computational time. For instance, a forward time integration of the 420-dimensional ODE approximation using \textsc{ode15s} of \textsc{matlab} (as depicted in Fig.~\ref{fig:mode1_timehistory_pre}) takes around 35 seconds, whereas the simulation with the two-dimensional SSM-based ROM~\eqref{eq:ROM1_stable} is completed in merely 1 second.

\subsubsection{Reduction for dynamics after bifurcation}
\label{subsec:5.1.2}
We set $\tau_d = 1.1>\tau^*$ and switch our focus to the dynamics after the Hopf bifurcation. In this case, the first five eigenvalues in ascending order of their real parts are updated as (cf.~\eqref{eq:lamd-duffing-before-hopf}).
\begin{equation}\lambda_{1,2} = 0.0102 \pm 1.5160\mathrm{i},\quad\lambda_3 = -2.5602 \pm 0.00\mathrm{i},\quad\lambda_{4,5} = -3.2606 \pm 6.6996\mathrm{i}.
\end{equation}
We take the unstable subspace as the master subspace for model reduction and compute the associated two-dimensional unstable SSM. The SSM-based ROM at $\mathcal{O}(9)$ is obtained below (cf.~\eqref{eq:red-auto}).
\begin{align}
\dot{\rho} &= -1.773\times10^{-10}\rho^9 -1.874\times10^{-8}\rho^7 -2.403\times10^{-6}\rho^5 -0.0004357\rho^3 +0.01023\rho, \nonumber\\
\dot{\theta} &= -1.212\times10^{-9}\rho^8 -1.34\times10^{-7}\rho^6 -1.78\times10^{-5}\rho^4 -0.003826\rho^2 +1.516.
\label{eq:ROM1_unstable}
\end{align}

\begin{figure}[ht] 
 \centering 
 \includegraphics[width=0.5\textwidth]{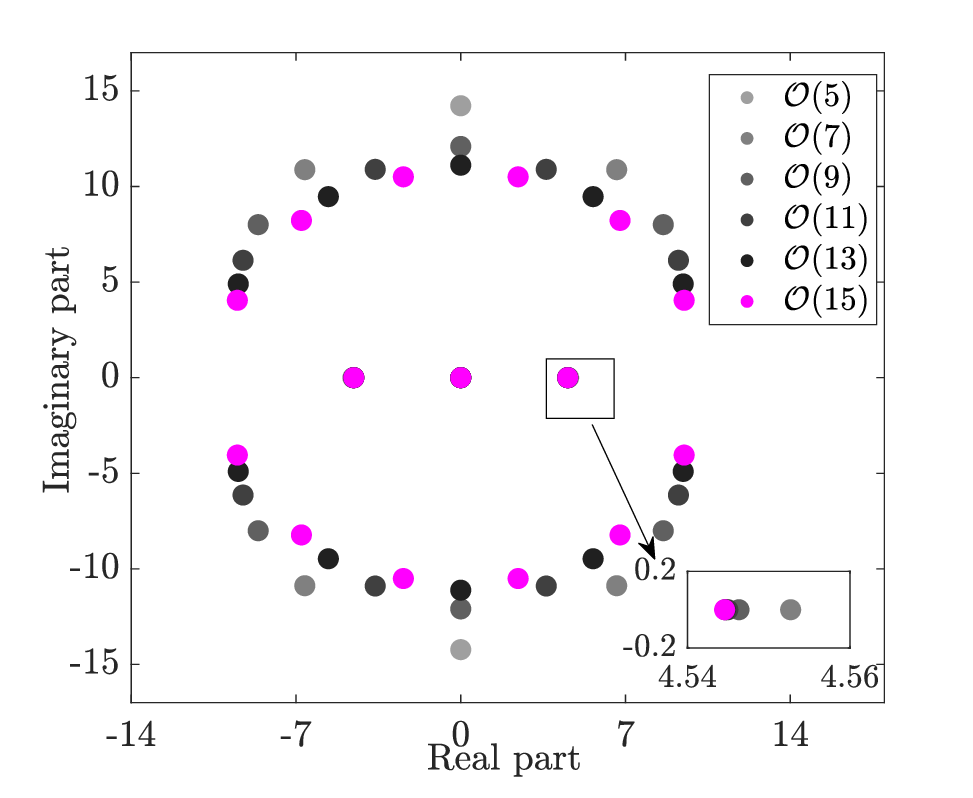} 
 \caption{\small Plot of the roots in the complex plane for $a(\rho)$ (cf. the first subequation of ~\eqref{eq:red-auto}) for the ROM~\eqref{eq:ROM1_unstable}, with darker colors indicating an increasing truncation order up to $\mathcal{O}(15)$. The zeros from the highest approximation are highlighted in magenta. We observe that a non-trivial transverse zero persists for higher-order truncation, and it is clearly within the domain of analyticity.} 
 \label{fig:model1_roots}
\end{figure}

From the polynomial function $a(\rho)=\dot{\rho}$ \textcolor{blue} { shown in the first subequation of~\eqref{eq:ROM1_unstable}}, it is known that there exists a non-trivial root $\rho^* \approx 4.54457$ such that $a(\rho^*) = 0$. This root $\rho^*$ corresponds to a limit cycle with frequency $b(\rho^*)$, as we have discussed in Sect.~\ref{sec:ssm-prediction}. We have calculated all the characteristic roots of $a(\rho)$ truncated at various expansion orders, as shown in Fig.~\ref{fig:model1_roots}. The darker colors represent the characteristic roots at higher orders, and the roots from the highest approximation are highlighted in magenta. We observe from Fig.~\ref{fig:model1_roots} that there is converged non-trivial root $\rho^* \approx 4.54457$.

\begin{figure}[!ht]
\centering
\includegraphics[width=.45\textwidth]{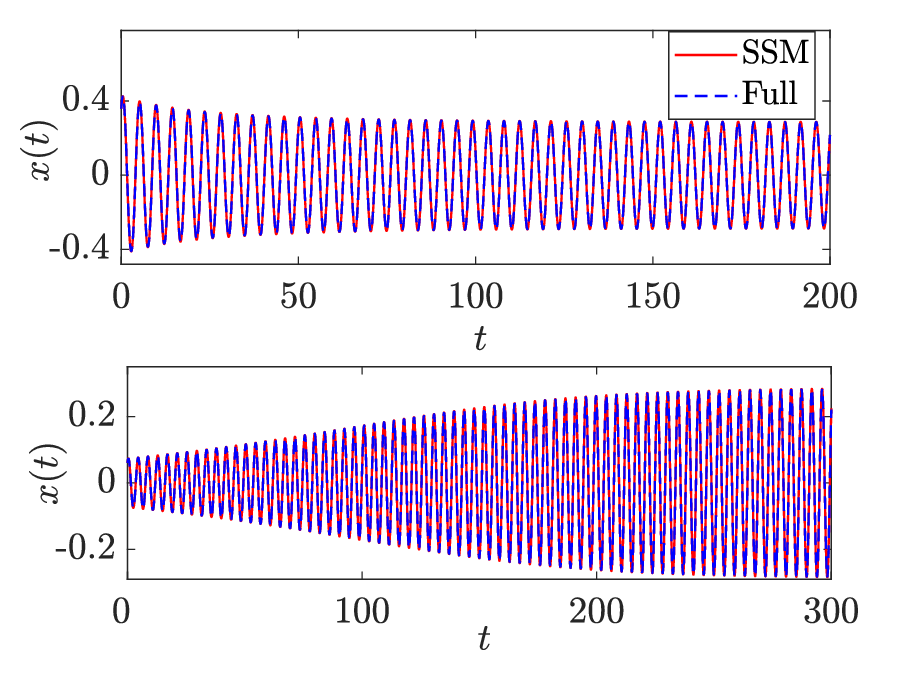}
\includegraphics[width=.45\textwidth]{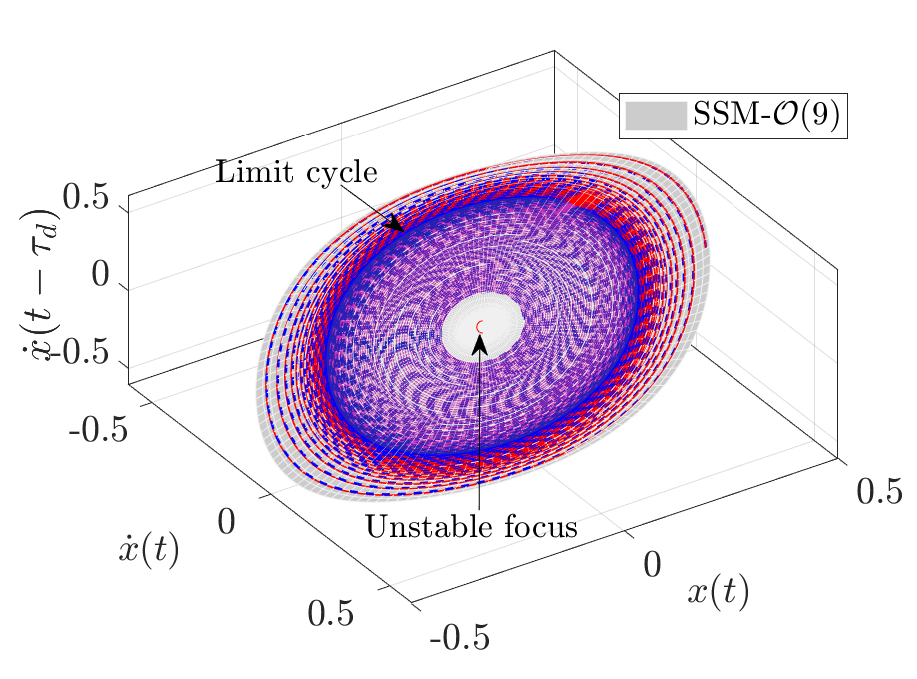}
\caption{\small Time history curve of the free variable $x(t)$ (left panel) and the projection of the SSM (gray plane) onto the coordinates ($\dot{x}(t), x(t), \dot{x}(t-\tau_d)$) (right panel). The blue dashed line and the red solid line correspond to the trajectories of the full system and that obtained via the ROM~\eqref{eq:ROM1_unstable}. The left panel demonstrates that the limit cycle attracts nearby trajectories.}
\label{fig:mode1_time_history_unstable}
\end{figure}

With the SSM computed up to $\mathcal{O}(9)$, the free vibration of the Duffing system was simulated under two different initial conditions. Initial condition 1 is given by $\bs{p}_0=(6.5 e^{0.1\mathrm{i}}, 6.5 e^{-0.1\mathrm{i}})$, with $\rho_0 = 6.5 > \rho^* = 4.54457$, and initial condition 2 is given by $\bs{p}_0=(1.2 e^{0.1\mathrm{i}}, 1.2 e^{-0.1\mathrm{i}})$, with $\rho_0 = 1.2 < \rho^*$. The free vibration response of the full system and that predicted via the ROM~\eqref{eq:ROM1_unstable} are shown in Fig.~\ref{fig:mode1_time_history_unstable}. From the first subplot in the left panel of Fig.~\ref{fig:mode1_time_history_unstable}, we see that the trajectory gradually approaches the limit cycle when the initial condition is outside the limit cycle. The second subplot shows that when the unstable focus is disturbed, the trajectory moves away from the focus and gradually approaches the limit cycle. Therefore, the SSM-based ROM performs excellently in predicting the response of the full system. In fact, the limit cycle stays on the two-dimensional SSM of the unstable focus, which is further verified in the right panel of Fig.~\ref{fig:mode1_time_history_unstable}.

\subsubsection{Forced vibration: isolas, bifurcations and quasi-periodic orbits}
\label{sec:5.1.3}
\begin{figure}[ht] 
 \centering 
\includegraphics[width=0.45\textwidth]{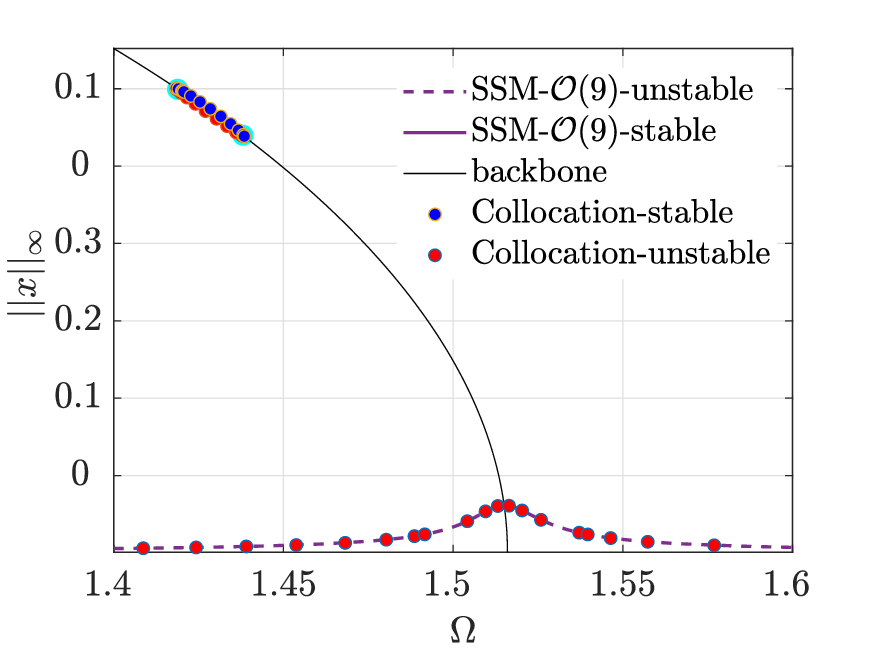}
\includegraphics[width=0.45\textwidth]{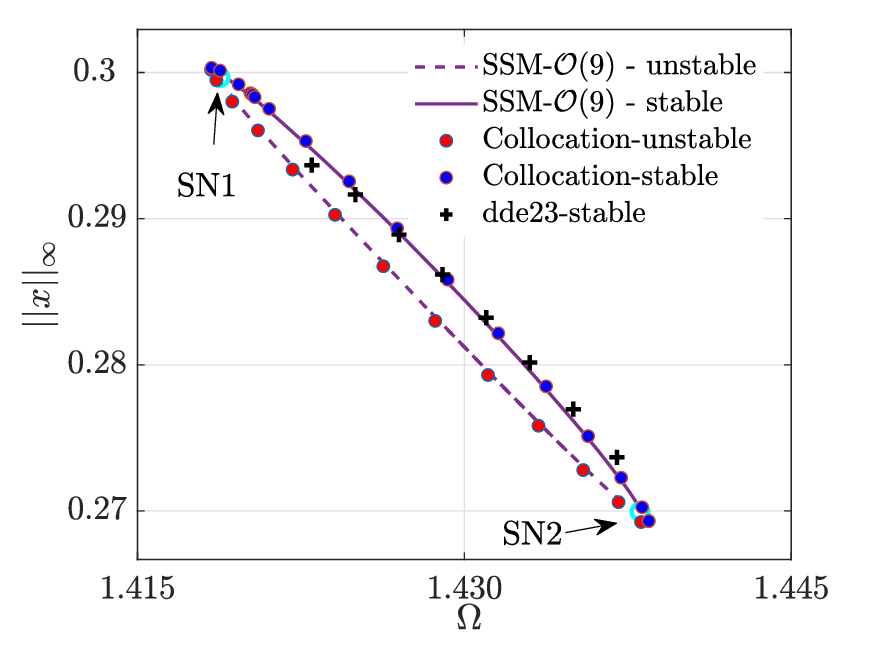} 
\caption{\small The forced response curves subjected to the periodic excitation term 0.0009 cos $\Omega t$ $(\epsilon = 0.0009)$ for the delayed Duffing system~\eqref{eq:dde-duffing}. The right panel is a zoomed plot of the left panel. The solid and dashed lines of the forced response curve of the periodic orbit represent stable and unstable periodic orbits, respectively, and the cyan circles indicate saddle-node bifurcation obtained from SSM-based prediction. Reference solutions are obtained from the collocation method applied to the ODE approximation of the DDEs. Results from the forward simulation of the DDEs are also provided in the right panel (dde23-stable).} 
\label{fig:model1_isola}
\end{figure}

We now move to the case of forced vibration, namely, $\epsilon>0$. In particular, we again take $\tau_d=1.1$ and focus on the dynamics after the bifurcation because the forced vibration in the post-bifurcation regime is more complicated than that of the pre-bifurcation regime. As illustrated in the previous subsection, $\rho^*$ is a nontrivial root of the polynomial $a(\rho)$, which corresponds to the limit cycle of the full system. This root corresponds to a point on the backbone curve of the unforced system ($\epsilon = 0$). Following~\cite{ponsioen2019analytic}, this point on the backbone curve is perturbed into an isola for $\epsilon > 0$ small enough. Indeed, as shown in Fig.~\ref{fig:model1_isola}, we obtain an isolated branch of periodic orbits with large amplitudes when $\epsilon=0.0009$. In particular, the isola has a segment of stable periodic orbits and a segment of unstable periodic orbits. These two segments are merged at two saddle-node (SN) points. In contrast, the periodic orbits on the whole main branch are unstable.

\begin{figure}[!ht]
 \centering 
\includegraphics[width=.45\textwidth]{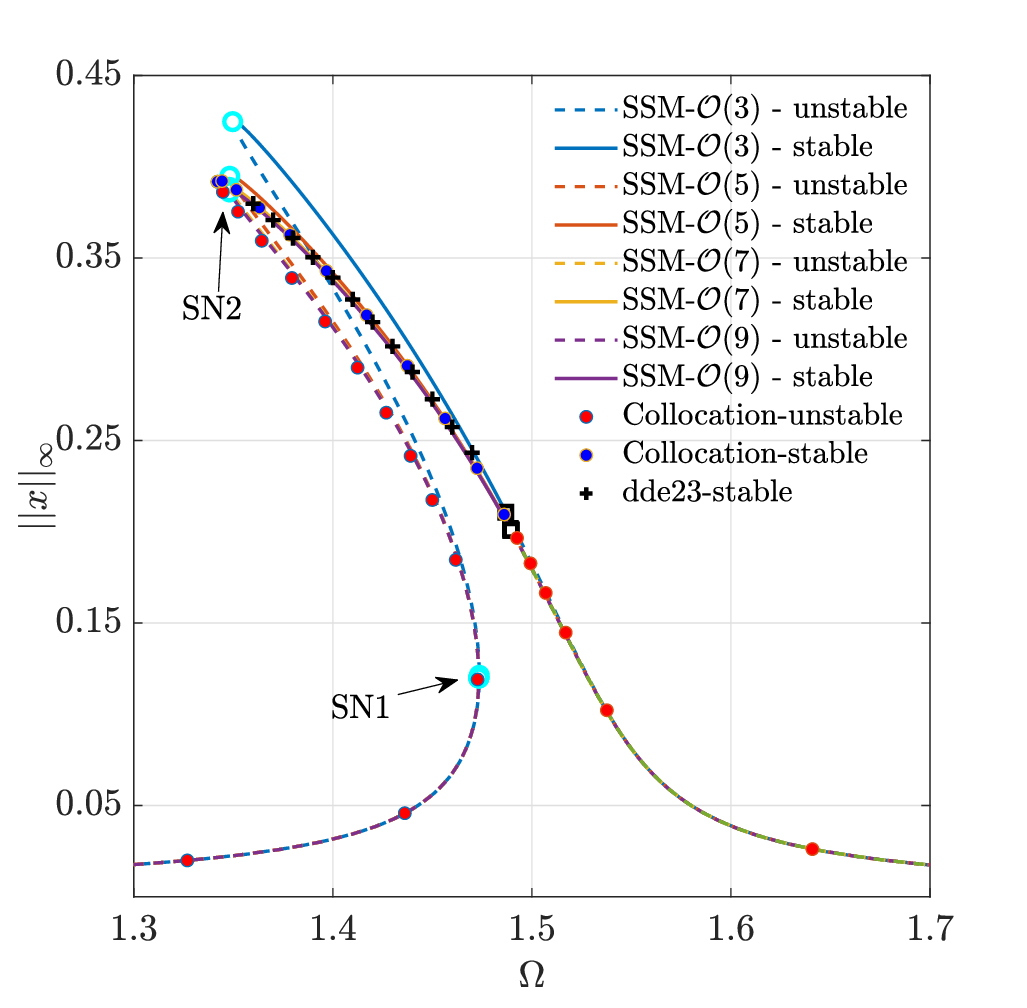} 
\includegraphics[width=.45\textwidth]{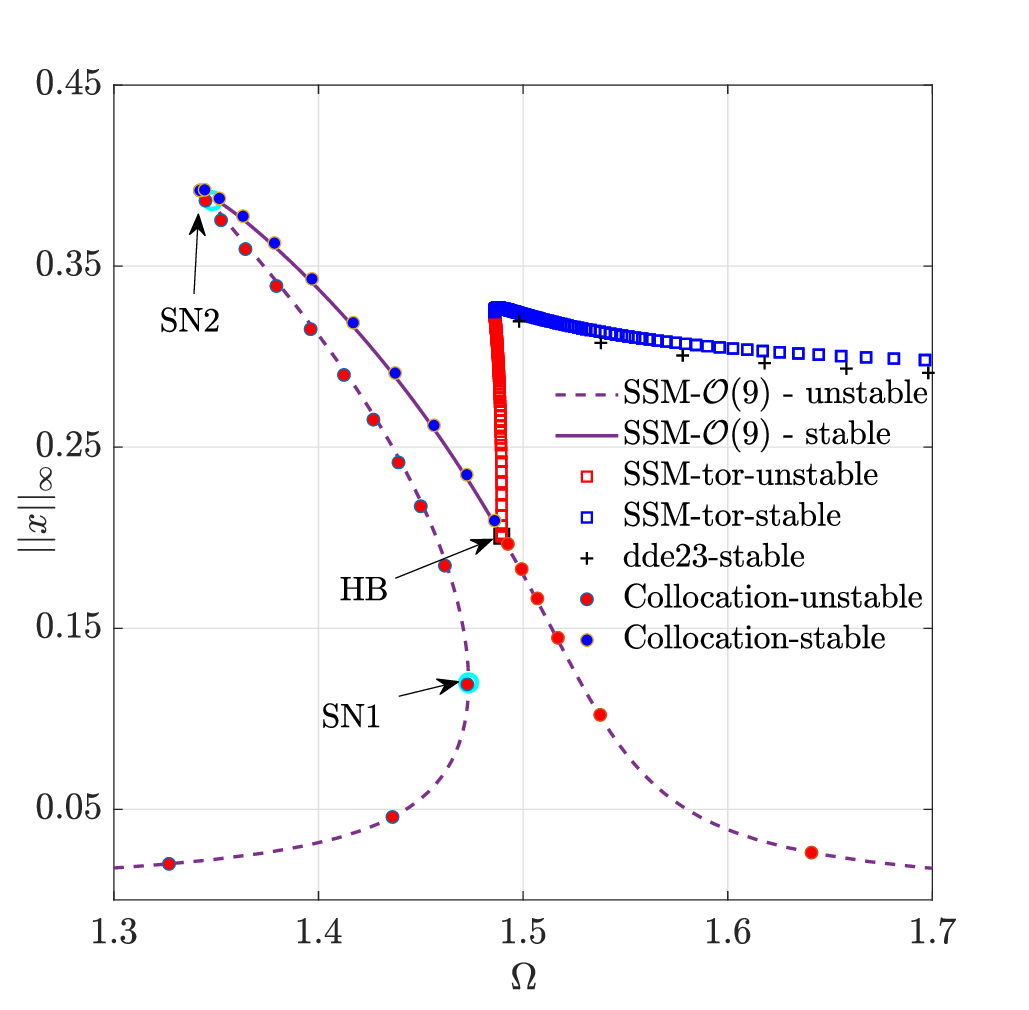} 
\caption{\small The forced response curves of the periodic orbits (left panel) along with quasi-periodic orbits (right panel) of the delayed Duffing system~\eqref{eq:dde-duffing} with $\epsilon=0.01$. The solid and dashed lines of the forced response curve of the periodic orbit represent stable and unstable periodic orbits, respectively, and the cyan circles indicate saddle-node bifurcation and the black squares indicate the Hopf bifurcations (HB). In the right panel, the blue and red squares denote stable and unstable quasi-periodic orbits predicted from SSM-based reduction, and the black crosses represent reference solutions obtained from forward simulation of the DDEs.}
\label{fig:demo1_FRC}
\end{figure}

To demonstrate the effectiveness and efficiency of SSM-based reduction, we calculate the forced response curve of the full, high-dimensional ODE approximation using the collocation method implemented in the po-toolbox of \textsc{coco}~\cite{dankowicz2013recipes,Schilder2024COCO,Ahsan2022COCO}. The obtained results are denoted as `Collocation' in Fig.~\ref{fig:demo1_FRC}. We further use a \textsc{matlab} solver `dde23' to perform forward simulation for the original DDEs~\eqref{eq:dde-duffing} to validate the predicted stable periodic orbits. As seen in Fig.~\ref{fig:demo1_FRC}, the results of the collocation method and forward simulation match well with those of the SSM-based reductions. Here, the computational time for the SSM-based reduction is about 10 seconds (for 185 points), while that of the collocation method is 51 hours and 49 minutes (yielding 233 points), indicating that the prediction via SSM reduction can achieve significant speed-up gains.

We note that the isola will merge with the main branch via a simple bifurcation as the forcing amplitude $\epsilon$ increases~\cite{li2023model}. Indeed, as seen in Fig.~\ref{fig:demo1_FRC}, the isola has merged with the main branch and disappeared when $\epsilon$ is increased to 0.01. Along this forced response curve (FRC), there are SN bifurcated periodic orbits, but also a secondary Hopf (or Neimark-Sacker) bifurcated periodic orbit. This Hopf bifurcation marks the birth of quasi-periodic orbits that stay on some tori. As demonstrated in~\cite{li2022nonlinear2}, such tori correspond to limit cycles of the SSM-based ROM~\eqref{eq:red-non-auto}. So, we further perform the continuation of limit cycles of the SSM-based ROM under variations in the excitation frequency, map these limit cycles to the tori of the full system, and extract the amplitude of the quasi-periodic orbits that stay on these tori. The obtained forced response curve regarding these quasi-periodic orbits is plotted along with the curve of a periodic orbit in the right panel of Fig.~\ref{fig:demo1_FRC}, where stable and unstable quasi-periodic orbits are denoted by red and blue squares. We observe the coexistence of stable quasi-periodic orbit and unstable periodic orbit for $\Omega>\Omega_\mathrm{HB}$. Thus, the forced vibration will approach a quasi-periodic orbit that stays on a torus in a steady state when $\Omega>\Omega_\mathrm{HB}$.

\begin{figure}[!ht]
\centering
\includegraphics[width=.45\textwidth]{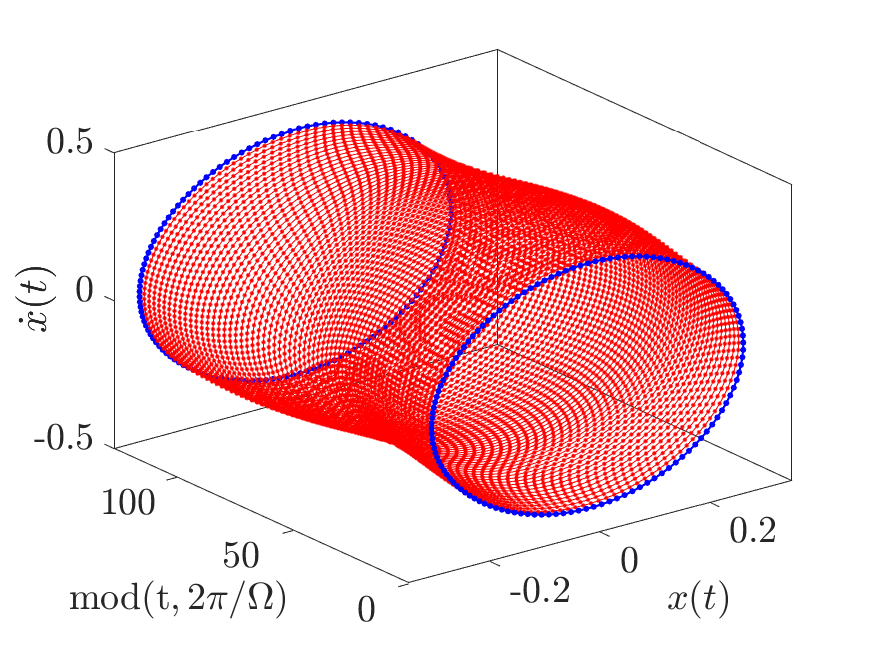}
\includegraphics[width=.45\textwidth]{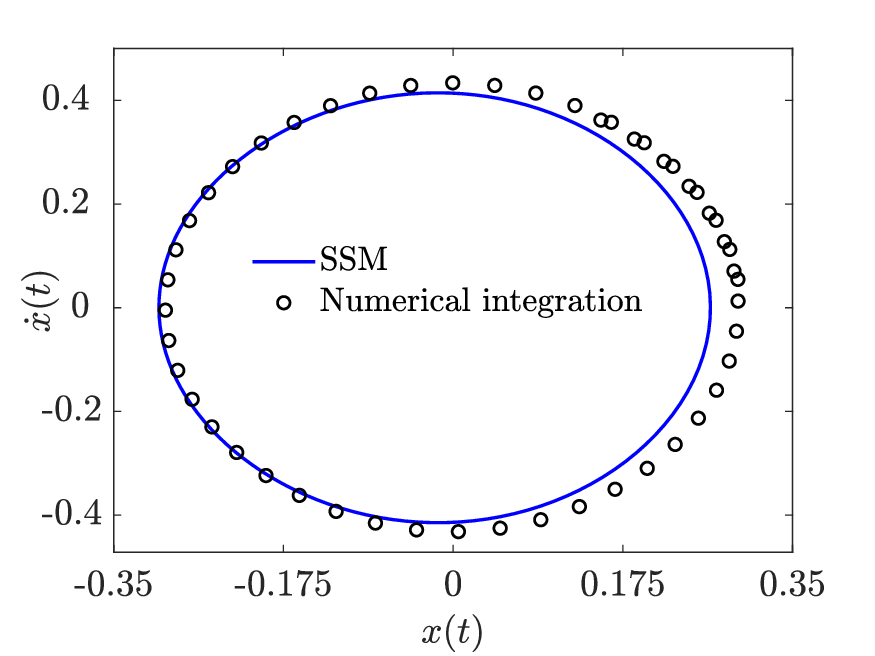}
\caption{\small (Left panel) A torus of the delayed Duffing system~\eqref{eq:dde-duffing} with $(\Omega,\epsilon) = (1.615, 0.01)$, obtained from SSM reduction (surface plot), and quasi-periodic orbits obtained from direct numerical integration of the ODE approximation system of~\eqref{eq:dde-duffing} shown as the red lines with dots. (Right panel) The intersection of the torus and the trajectory with a period-$2\pi/\Omega$ Poincaré section.}
\label{fig:mode1_Tor}
\end{figure}

To validate the effectiveness of the SSM-based reductions, we again apply the collocation method to the approximated ODE system to extract reference solutions for the forced periodic orbits. We also apply dde23 to perform forward simulations of the DDEs~\eqref{eq:dde-duffing} for some sampled $\Omega$ within the interval $[\Omega_\mathrm{SN1},\Omega_\mathrm{HB}]$. As we can see from the left panel of Fig.~\ref{fig:demo1_FRC}, the $\mathcal{O}(9)$ SSM-based predictions for periodic orbits match well with the reference solutions from the collocation and forward simulation methods. We further apply the forward simulation to the DDEs~\eqref{eq:dde-duffing} to extract the stable quasi-periodic orbits in a steady state. We plot the FRC for periodic and quasi-periodic orbits in the right panel of Fig.~\ref{fig:demo1_FRC}. A sampled torus, along with the quasi-periodic orbit in a steady state, is presented in Fig.~\ref{fig:mode1_Tor}. We observe that the SSM-based predictions for the quasi-periodic orbits also match the reference solutions. Here, the averaged computational time for numerical integration to yield a quasi-periodic orbit is about 10.85 seconds (refer to the '+' shaped marks in Fig.~\ref{fig:demo1_FRC}). In contrast, the averaged computational time for a torus via the SSM-based prediction is about 0.25 seconds per point. This again demonstrates the significant speed-up gain of SSM-based model reduction.

\subsubsection{Limitations of SSM-based reduction}
\label{sec:limitation}
We conclude this example with illustrations of the limitations of our SSM-based reduction. In particular, as we discussed in Remark \ref{rmk:ssm-reduction}, the response must be within the domain of convergence of Taylor expansion of SSM parametrization. This indicates that the SSM-based reduction may fail to capture the limit cycle in the post-flutter region if $\tau_d$ is much larger than $\tau^\ast=1.035$, where the Hopf bifurcation occurs.

To illustrate the aforementioned failure, we take $\tau_d=1.75$, keep all other parameters unchanged, and repeat the analysis in Section \ref{subsec:5.1.2}. The obtained results are detailed in ~\ref{appendix:C}. Indeed, as seen in Fig.~\ref{fig:root_plot_fail}, the SSM-based ROM fails to capture a non-spurious, non-trivial root of $a(\rho)=0$, and the non-trivial root $\rho_1^+$ is on the boundary of the domain of convergence. As expected, this spurious root $\rho_1^+$ makes an incorrect prediction to the limit cycle~\cite{ponsioen2019analytic}, as seen in Fig.~\ref{fig:root_plot_fail}.

We further illustrate the limitation of our SSM-based reduction for the case of non-autonomous dynamics. We revisit Section \ref{sec:5.1.3} but increase the forcing amplitude $\epsilon$. We have chosen $\epsilon=0.009$ and $\epsilon=0.01$ in Section \ref{sec:5.1.3}, and the SSM-based prediction yields convergent FRCs with increasing expansion orders for these two forcing levels, as seen in Figures~\ref{fig:model1_isola} and~\ref{fig:demo1_FRC}. However, the SSM-based reduction produces divergent forced periodic orbits with large amplitudes when $\epsilon=0.4$, as seen in Fig.~\ref{fig:frc-plot-fail}. This failure has also been observed and discussed in mechanical systems with internal resonance~\cite{li2022nonlinear1}.

\subsection{Two coupled oscillators with delay}
\label{sec:ex-two-os}
As our second example, we consider two coupled nonlinear oscillators with delay~\cite{saeed2017utilizing}. The system is subject to a base excitation, and the delay is used to quench the vibration of the system. The equations of motion of this system are listed below.

\begin{align}
\label{eq:model2_two_order}
\ddot{u}+&\mu_1\dot{u}+\omega_1^2u+2\gamma v u+\gamma v^2u+\gamma u^3 = \epsilon\Omega^2\cos\Omega t-(\beta_1u(t-\tau_d)+\beta_2\dot{u}(t-\tau_d)),\nonumber\\
\ddot{v}+&\mu_2\dot{v}+\omega_2^2v+\gamma v^2+3\gamma v^2+\gamma u^2v+\gamma v^3 = \epsilon\Omega^2\sin\Omega t-\left(\beta_1v(t-\tau_d)+\beta_2\dot{v}(t-\tau_d)\right).
\end{align}
Let ${q_1}=u,{q_2}=\dot{u},{q_3}=v$, and ${q_4}=\dot{v}$,~\eqref{eq:model2_two_order} can be rewritten in the first-order form~\eqref{eq: DDEs_form} below
\begin{align}
\dot{q_1} &= q_2,\nonumber\\
\dot{q_2} &= -\mu_1 q_2 - \omega_1^2 q_1 - 2\gamma q_1 q_3 - \gamma q_1 q_3^2 - \gamma q_1^3 + \epsilon \Omega^2 \cos\Omega t - (\beta_1 q_1(t-\tau_d) + \beta_2 q_2(t-\tau_d)), \nonumber\\
\dot{q_3} &= q_4, \nonumber\\
\dot{q_4} &= -\mu_2 q_4 - \omega_2^2 q_3 - \gamma q_3^2 - 3\gamma q_3^2 - \gamma q_1^2 q_3 - \gamma q_3^3 + \epsilon \Omega^2 \sin\Omega t - (\beta_1 q_3(t-\tau_d) + \beta_2 q_4(t-\tau_d)).
\label{eq:model2_one_order}
\end{align}

In the following computations, parameters are set as $\mu_1 = 0.015$, $\mu_2 = 0.035$, $\omega_1 = \sqrt{1+\gamma}$, $\omega_2 = \sqrt{1+3\gamma}$, $\gamma=0.3$, $\beta_1=-0.3$, $\beta_2=-0.1$, $\tau_d = 0.5$, and $\epsilon=0$ unless otherwise stated. Likewise, we follow Section~\ref{sec:dde-to-ode} to transform the original DDEs~\eqref{eq:model2_one_order} into a system of ODEs~\eqref{eq:non_eq1}. Here, we have $n=4$, and hence $\bs{z}\in\mathbb{R}^{4(2N+1)}$. As detailed in~\ref{sec:app-ex2}, we can choose $N=20$ to yield converged solutions. Thus, we take $N=20$ in this example. Moreover, the analysis in~\ref{sec:app-ex2} shows that the system undergoes a Hopf bifurcation at $\beta_1=\beta_1^\ast=-0.146$, as seen in the right panel of Fig.~\ref{fig:model2_r_distributed}.

\subsubsection{Reduction for dynamics before bifurcation}
By setting $\beta_1=-0.3<\beta^*$, we study the pre-bifurcation dynamics of system~\eqref{eq:model2_two_order}. The first three pairs of eigenvalues of the linear portion of the approximated system, corresponding to Eq.~\eqref{eq:model2_two_order}, are:
 The first 8 nonzero eigenvalues are given as 

\begin{align}
\lambda_{1,2} = -0.0351 \pm 0.9936\mathrm{i}, \quad
\lambda_{3,4} = -0.0468 \pm 1.2571\mathrm{i}, \quad
\lambda_{5,6} = -9.9286 \pm 7.4053\mathrm{i}. \quad
\label{eq:lamd-model2-before-hopf}
\end{align}

We compute the SSM corresponding to the first pair of eigenvalues and formulate an ROM on the SSM. In the scenario of free vibration where $\epsilon=0$, we derive a two-dimensional ROM~\eqref{eq:ROM2_stable} for the SSM approximated at $\mathcal{O}(9)$. The ROM is listed as follows (cf.~\eqref{eq:red-auto}):

\begin{align}
\label{eq:ROM2_stable}
\dot{\rho} &=-2.133 \times 10^{-10}  {\rho}^9 + 9.021 \times 10^{-9}  {\rho}^7 - 2.574 \times 10^{-7}  {\rho}^5 - 3.25 \times 10^{-5}  {\rho}^3 - 0.03514  \rho, \nonumber\\
\dot{\theta} &=-7.883 \times 10^{-9}  {\rho}^8 + 5.242 \times 10^{-7}  {\rho}^6 - 4.172 \times 10^{-5}  {\rho}^4 + 0.005429  {\rho}^2 + 0.9936.
\end{align}

\begin{figure}[!ht]
\centering
\includegraphics[width=.45\textwidth]{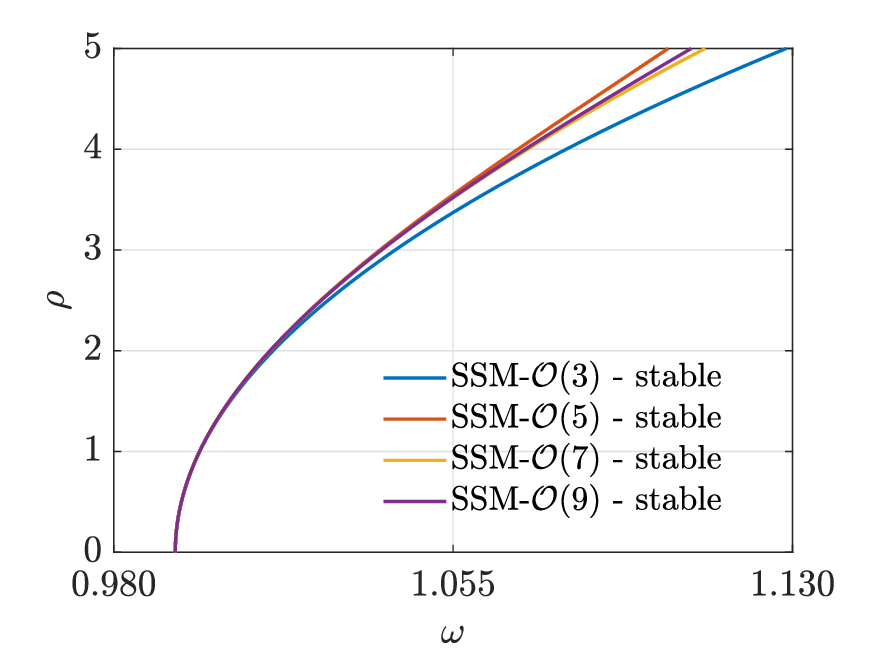} 
\includegraphics[width=.45\textwidth]{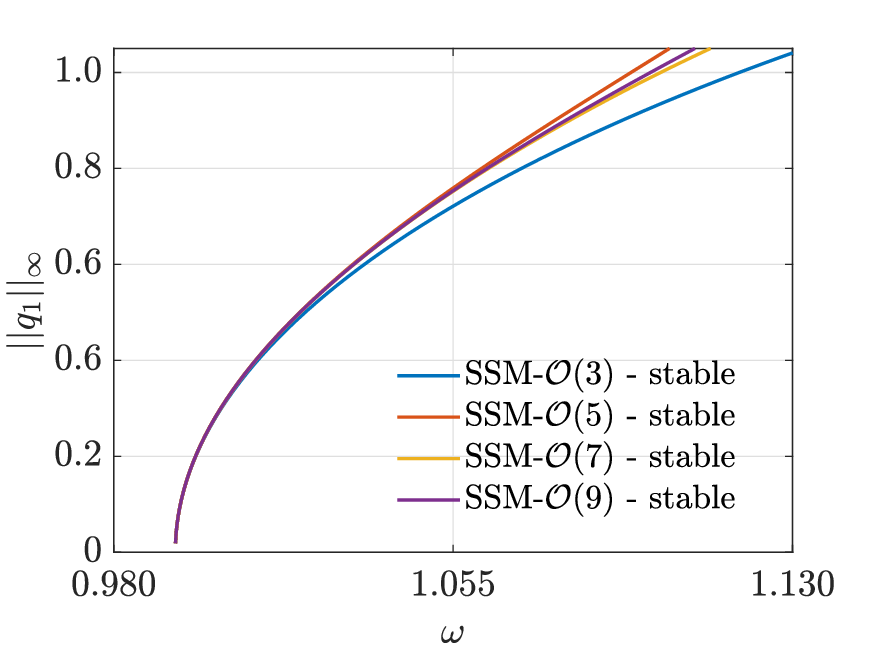}
\caption{\small Backbone curves in polar (reduced) and physical coordinates under increasing orders of approximation for the slowest two-dimensional SSM of the system of two coupled oscillators~\eqref{eq:model2_one_order}.}
\label{fig:model2_backbone}
\end{figure}

Similarly to the previous examples, we obtain damped backbone curves in polar (reduced) and physical coordinates shown in Fig.~\ref{fig:model2_backbone}. We observe that higher-order expansions are needed to obtain convergence in backbones at higher amplitudes. In particular, the backbone curves in Fig.~\ref{fig:model2_backbone} show convergence at $\mathcal{O}(5)$ for $\rho \leq 2.32$, at $\mathcal{O}(7)$ expansion for $\rho \leq 3.52$, and at $\mathcal{O}(9)$ expansion for $\rho \leq 4.52$, Here, the backbone curves exhibit hardening behavior.

\begin{figure}[!ht]
\centering
\includegraphics[width=.45\textwidth]{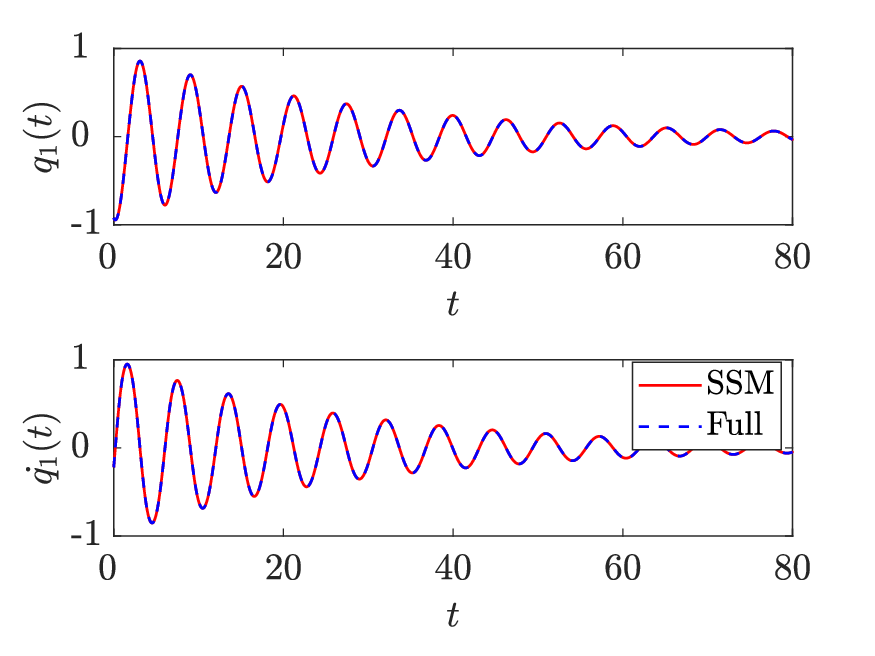} 
\includegraphics[width=.45\textwidth]{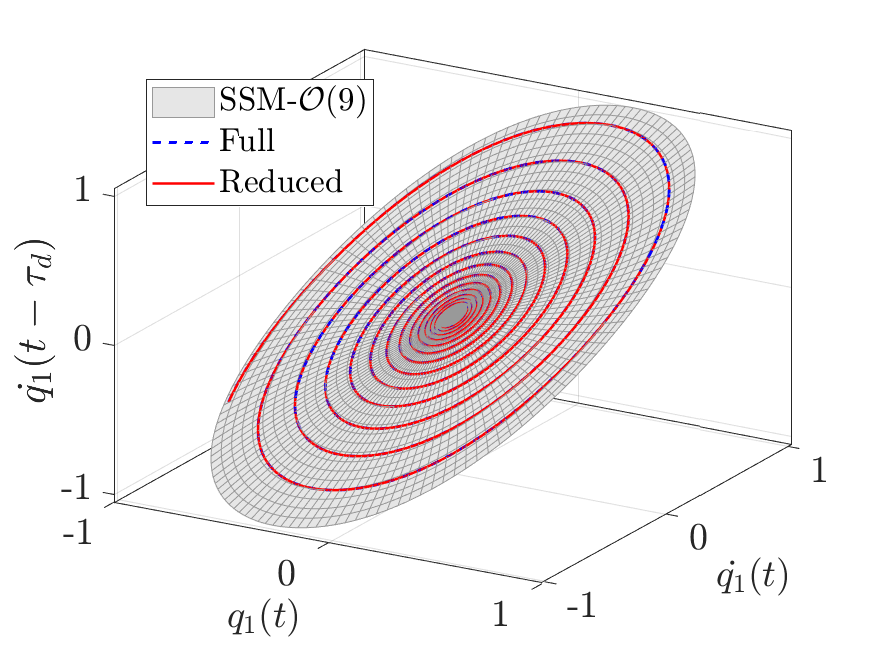}
\caption{\small (Left panel) Time history plots of $q_1(t)$ and $\dot{q_1}(t)$ for system~\eqref{eq:model2_two_order}. (Right panel) Visualization of the SSM projected onto the coordinates ($q_1(t), \dot{q_1}(t),\dot{q_1}(t-\tau_d)$) (gray plane) along with the trajectories of the full system (blue dashed line) and the ROM~\eqref{eq:ROM2_stable} (red solid line).}
\label{fig:model2_time_history}
\end{figure}

Now, we use the SSM-based ROM to predict the free vibration of the system. We select an initial state on the SSM and perform forward simulations for both the full system~\eqref{eq:non_eq1} and the SSM-based ROM~\eqref{eq:ROM2_stable}. The initial conditions on the SSM are set as $\bs{p}_0=(4.5 e^{0.1\mathrm{i}}, 4.5 e^{-0.1\mathrm{i}})$. The simulated results, depicted in Fig.~\ref{fig:model2_time_history}, demonstrate that the trajectory from the full system simulation aligns with the trajectory predicted from the SSM-based ROM~\eqref{eq:ROM2_stable}. Furthermore, we project the SSM onto the coordinates ($q_1(t), \dot{q_1}(t), \dot{q_1}(t-\tau_d)$), as illustrated in Fig.~\ref{fig:model2_time_history}. The full-system trajectory stays on the SSM and overlaps with the prediction of the SSM-based ROM. This confirms the invariance of the computed SSM.

\begin{figure}[!ht]
 \centering 
 \includegraphics[width=.5\textwidth]{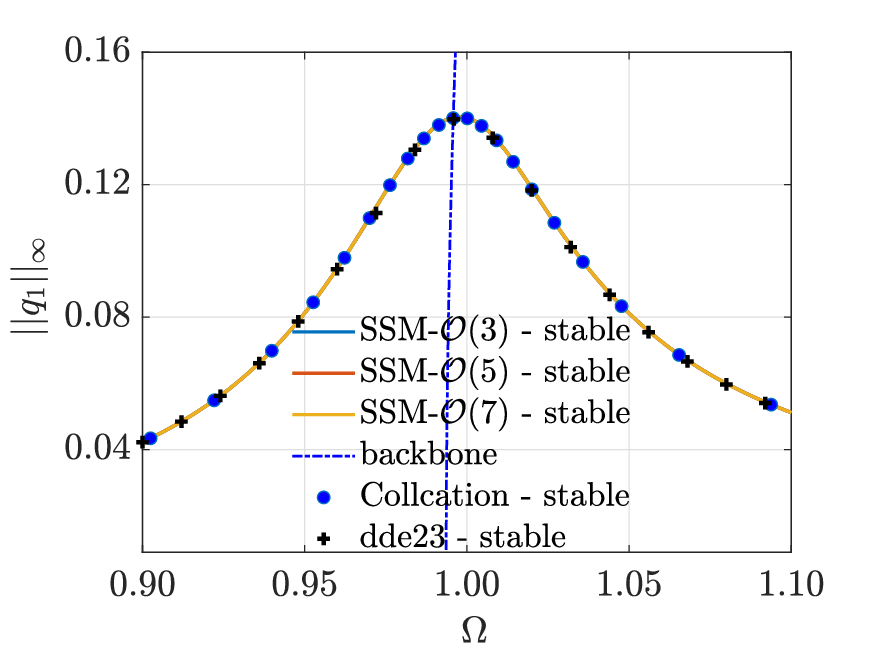} 
 \caption{\small The forced response curve of periodic orbits of the two coupled oscillators~\eqref{eq:model2_two_order} with $\epsilon = 0.01$.} 
 \label{fig:model2_FRCs}
\end{figure}

We further investigate the forced vibration of the system with $\epsilon=0.01$, with particular emphasis on the system's forced response curve (FRC) near its first natural frequency (i.e., $\Omega \approx \omega_1$). The FRC under different expansion orders is shown in Fig.~\ref{fig:model2_FRCs}, where it is observed that $\mathcal{O}(3)$ expansion is already sufficient. To validate the effectiveness of the SSM-based reductions, we applied the collocation method to the approximated ODE system of~\eqref{eq:model2_one_order}, extracting reference solutions for the forced periodic orbits. Additionally, we utilized the \textsc{matlab} solver 'dde23' to conduct forward simulations of the original DDEs~\eqref{eq:model2_two_order} to validate the predicted stable periodic orbits. As shown in Fig.~\ref{fig:model2_FRCs}, the collocation method and forward simulation results align well with the SSM-based predictions.

\subsubsection{Reduction for dynamics after bifurcation}
We set $\beta_1=-0.145>\beta_1^*$ and focus on the post-bifurcation dynamics. In this context, the first six eigenvalues are obtained and arranged in ascending order of their real parts, as shown below
\begin{align}
\label{eq:lamd-model2-after-hopf}
\lambda_{1,2} = 0.0010 \pm 1.0590\mathrm{i}, \quad
\lambda_{3,4} = -0.0112 \pm 1.3058\mathrm{i}, \quad
\lambda_{5,6} = -9.6918 \pm 7.6345\mathrm{i}. \quad
\end{align}
We take the unstable spectral subspace corresponding to $\lambda_{1,2}$ as the master subspace of the SSM and compute the associated two-dimensional SSM. The SSM-based ROM at $\mathcal{O}(9)$ expansion is obtained below (cf.~\eqref{eq:red-auto}):
\begin{align}
\label{eq:ROM2_unstable}
\dot{\rho} &= 1.593 \times 10^{-11}\,\rho^9 - 1.249 \times 10^{-9}\,\rho^7 + 1.386 \times 10^{-7}\,\rho^5 - 4.157 \times 10^{-5}\,\rho^3 + 0.001039\,\rho, \nonumber\\
\dot{\theta} &=- 3.079 \times 10^{-9}\,\rho^8 + 2.616 \times 10^{-7}\,\rho^6 - 2.674 \times 10^{-5}\,\rho^4 + 0.004489\,\rho^2 + 1.059.
\end{align}

\begin{figure}[!ht]
\centering
\includegraphics[width=.45\textwidth]{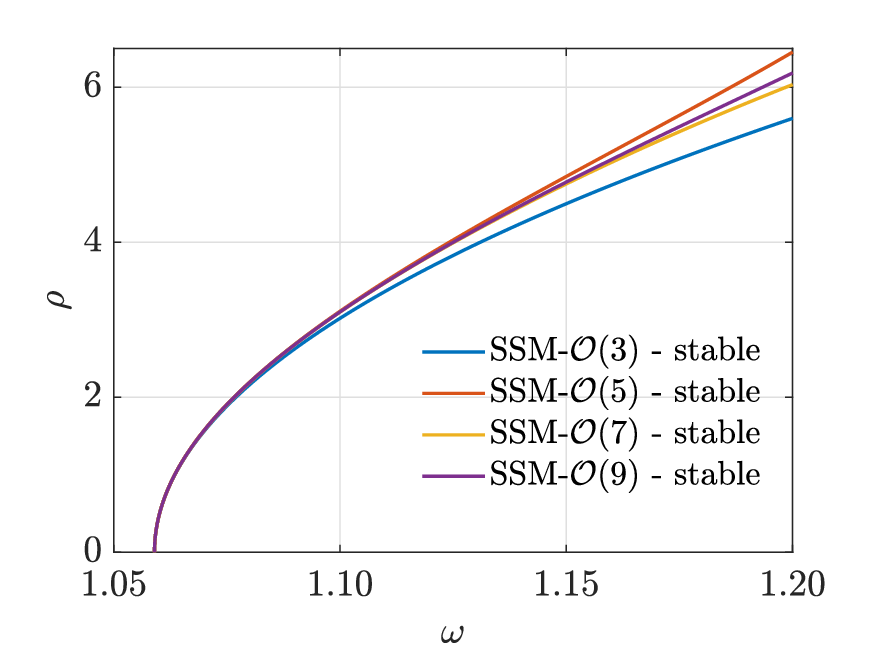} 
\includegraphics[width=.45\textwidth]{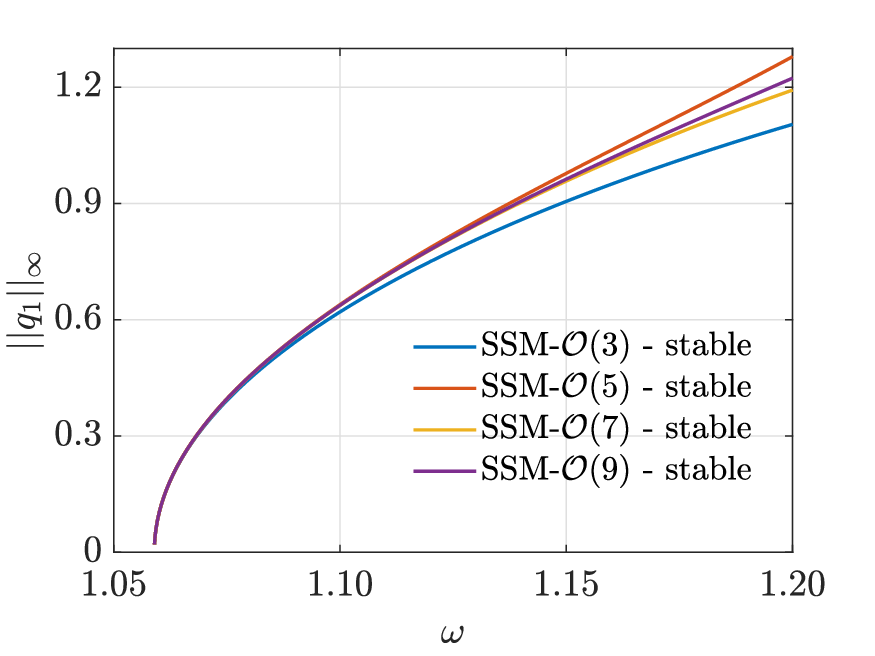}
\caption{\small Backbone curves in polar (reduced) and physical coordinates under increasing orders of approximation for the two-dimensional SSM of the system~\eqref{eq:ROM2_unstable}.}
\label{fig:model2_backbone_post}
\end{figure}

We first consider the free vibration of the two coupled oscillators~\eqref {eq:model2_two_order}, i.e., $\epsilon=0$. Similar to the previous examples, we have computed the damped backbone curve via two-dimensional SSM of system~\eqref{eq:ROM2_unstable}. Fig.~\ref{fig:model2_backbone_post} shows the backbone curve of ROM~\eqref{eq:ROM2_unstable} at various expansion orders of SSM approximation. This backbone curve again displays hardening behavior. As shown in the left panel of Fig.~\ref{fig:model2_backbone_post}, the backbone curve is well converged with an expansion up to $\mathcal{O}(9)$ for amplitude $\rho \leq 5.4$.

\begin{figure}[!ht]
\centering
\includegraphics[width=0.9 \textwidth]{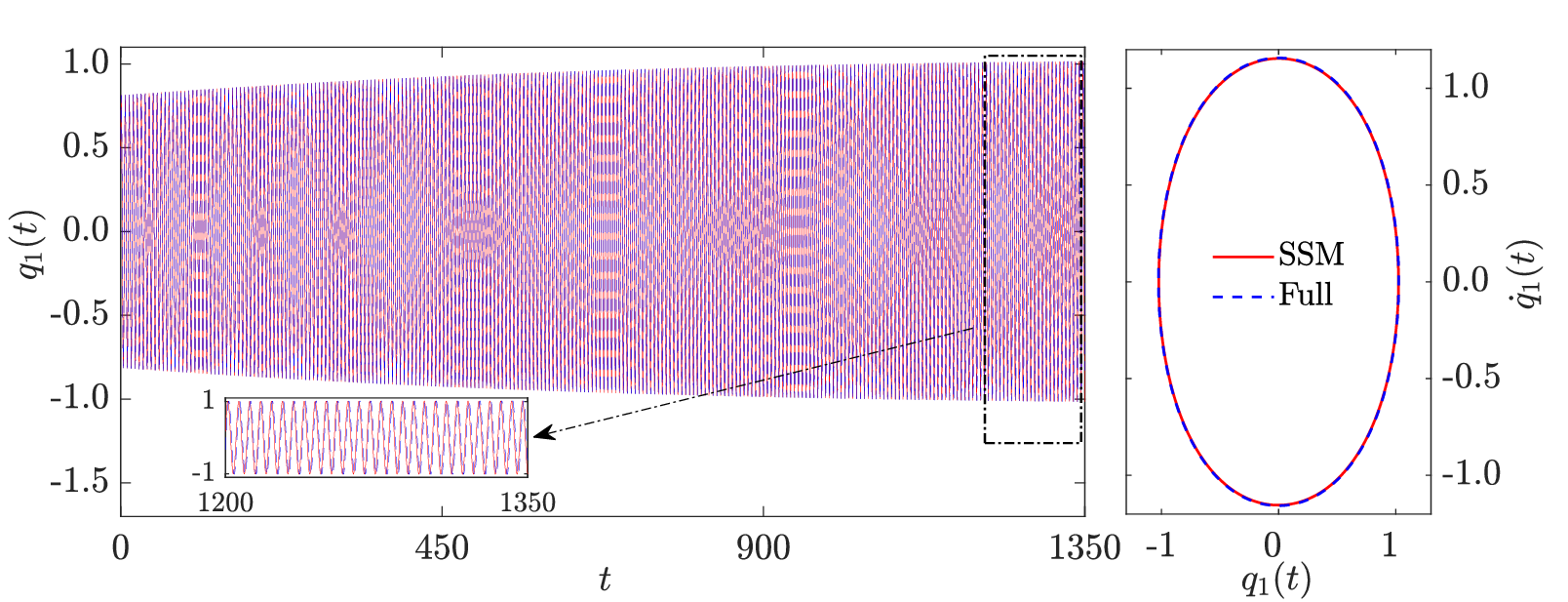}
\caption{\small Time history curve of the free variable $q_1(t)$. The projections of the trajectories of the full system (blue dashed line) and the ROM~\eqref{eq:ROM2_unstable} (red solid line)are depicted. In the left subplot, the trajectories approach the limit cycle in the steady state, while the right subplot illustrates the limit cycle of the system in its steady state. }
\label{fig:model2_limitcycle}
\end{figure}

We observe from the SSM-based ROM~\eqref{eq:ROM2_unstable} that the polynomial function $a(\rho)$ (cf.~\eqref{eq:red-auto}) admits a nontrivial root $\rho^*=5.179$ such that $a(\rho^*)=0$. This root $\rho^*$ corresponds to a limit cycle with frequency $b(\rho^*)$, as discussed in Sect.~\ref{sec:ssm-prediction}. We take the initial condition $\bs{p}_0=(4.0 e^{0.1\mathrm{i}}, 4.0 e^{-0.1\mathrm{i}})$, and perform forward time integration for both the ROM~\eqref{eq:ROM2_unstable} and the approximate system~\eqref{eq:model2_one_order} using the same initial conditions. As shown in Fig.~\ref{fig:model2_limitcycle}, the trajectories indeed approach the limit cycle in the steady state (as depicted in the left subplot). Furthermore, we plot the limit cycle in the right subplot. The SSM-based ROM provides excellent predictions, with the full system's simulated trajectory aligning with the ROM's trajectory, confirming the accuracy of the SSM-based predictions.

\begin{figure}[!ht]
\centering
\includegraphics[width=.45\textwidth]{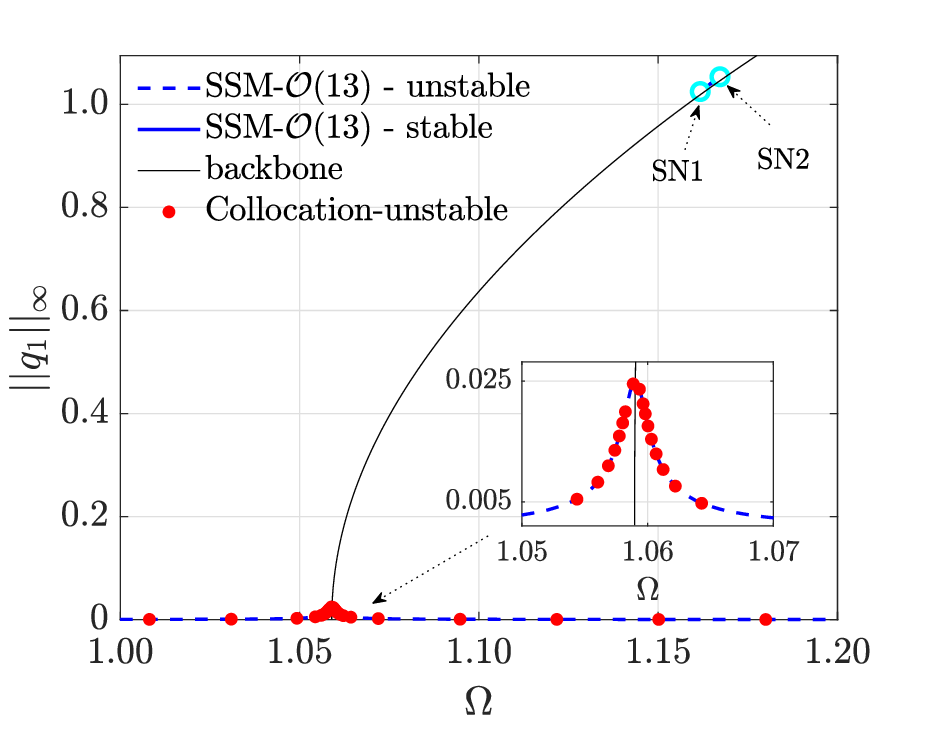}
\includegraphics[width=.45\textwidth]{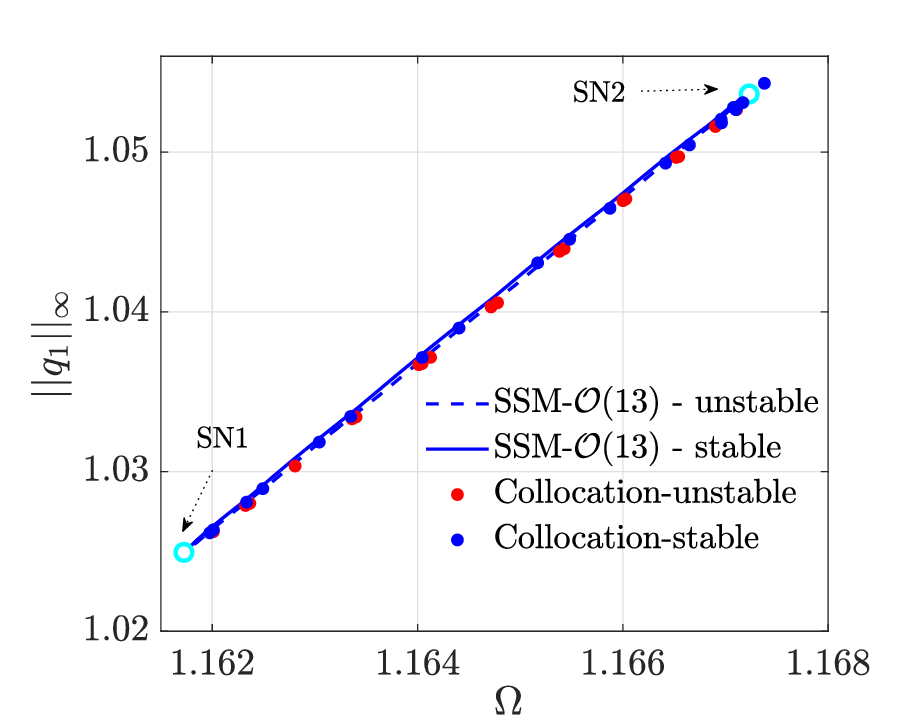}
\caption{\small The forced response curves subjected to the periodic excitation term $\epsilon$ cos $\Omega t$ for the two coupled oscillator~\eqref{eq:model2_one_order}. (Left panel) For $\epsilon = 0.00005$, the system exhibits a small amplitude response on the main branch and a large amplitude response on an isolated branch. A enlarged plot for the isola is given in the right panel. The solid and dashed lines represent stable and unstable periodic orbits, respectively. Cyan circles indicate saddle-node bifurcations predicted using SSM-based methods. Reference solutions are derived from the collocation method applied to the ODE approximation system of the DDEs.}
\label{fig:model2_FRC_isola}
\end{figure}

Next, we analyze the forced vibration of the system, i.e., when $\epsilon > 0$. In this poster-flutter regime, the dynamic characteristics of the forced vibration become more complex than before. As mentioned in previous sections, $\rho^*$ is a nontrivial solution of the polynomial $b(\rho)$, corresponding to a limit cycle of the system located on the backbone curve of the unforced system ($\epsilon = 0$). Following the findings of~\cite{ponsioen2019analytic}, the point on the backbone curve is perturbed to form an isola when $\epsilon>0$ is small enough. Specifically, as shown in the left panel of Fig.~\ref{fig:model2_FRC_isola}, when $\epsilon = 0.00005$, we observe a small-amplitude response on the main branch and a larger-amplitude response on an isolated branch. The isolated branch contains a segment of stable periodic orbit and a segment of unstable periodic orbit, which intersect at two saddle-node (SN) points. In contrast, the periodic orbits on the main branch remain unstable.

Once again, we apply the po toolbox in \textsc{coco} to extract the FRC of the full nonlinear system as the reference solution to compare the accuracy of solutions obtained by SSM reduction methods. The results from po are labeled as Collocation. Fig.~\ref{fig:model2_FRC_isola} shows that the results from the collocation method align closely with those from the SSM-based reduction. The computation time for the SSM-based reduction was approximately 10 seconds, whereas the collocation method took 38 hours and 41 minutes, demonstrating a significant acceleration achieved through SSM-based reduction.

\subsection{Hutchinson equation}
\label{sec:ex-hutchinson}
In our final example, we consider the Hutchinson equation with a diffusion term~\cite{faria2000normal}. The equation is in the form of a delay partial differential equation (DPDE) and is given below:
\begin{equation} 
\label{eq:model3_PDE}
\dfrac{\partial u(t,x)}{\partial t}=d\dfrac{\partial^2u(t,x)}{\partial x^2}+u(t,x)-au(t-1,x)-au(t-1,x)u(t,x),
\end{equation}
with boundary condition
\begin{equation}
 u(t,0)=u(t,\pi)=0.
\end{equation}
Here, $u(t,x)$ represents the population size at time $t$ and position $x$. This equation models the population dynamics incorporating both diffusion and delay effects. Following~\cite{faria2000normal}, we apply the Galerkin method to discretize the DPDE in spatial domain, such that the DPDE is transformed into DDEs. In particular, we substitute the expansion below into~\eqref{eq:model3_PDE}
\begin{equation}
 u(t,x)=\sum_{k=1}^M \beta_k(x)q_k(t), \quad \beta_k(x)=\sqrt{\frac{2}{\pi}}\sin kx.
\end{equation}
and, apply a Galerkin projection, which yields
\begin{equation} 
\label{eq:model3_DDE_N}
\dot{q}_i(t)=-d\sum_{k=1}^Mq_k(t)\Big\langle\beta_i(x),k^2\beta_k(x)\Big\rangle + q_i(t) -aq_i(t-1)-a\sum_{j=1}^M\sum_{k=1}^Mq_j(t-1)q_k(t)\Big\langle\beta_i(x),\beta_j(x)\beta_k(x)\Big\rangle
\end{equation}
for $i=1,\cdots,M$. Here, $\langle\cdot, \cdot \rangle$ denotes the inner product as below
\begin{equation}
\label{eq:inner_product}
\langle f, g \rangle = \int_0^\pi f(x) g(x) \, dx.
\end{equation}
A spectral analysis is performed at the equilibrium point $\bs{q}=0$ for~\eqref{eq:model3_DDE_N}, which leads to the characteristic equation below:
\begin{equation}
\label{eq:model3_lambda_eq}
\lambda+ae^{-\lambda}+dk^2-1=0,\quad k=0,1,\cdots.
\end{equation}
In the rest computations, we choose $d = 1$ following~\cite{faria2000normal}. When $a < {\pi}/{2}$, all eigenvalues have real parts less than zero, indicating that the equilibrium point is stable. At $a = {\pi}/{2}$, the eigenvalues present a pair of purely imaginary roots, while the real parts of all other eigenvalues remain below zero, marking the onset of a Hopf bifurcation. Therefore, $a = {\pi}/{2} = a^*$ is identified as the bifurcation point. In fact, this is a supercritical Hopf bifurcation point~\cite{faria2000normal}.

\subsubsection{Reduction for dynamics before bifurcation}
\begin{figure}[!ht]
\centering
\includegraphics[width=.45\textwidth]{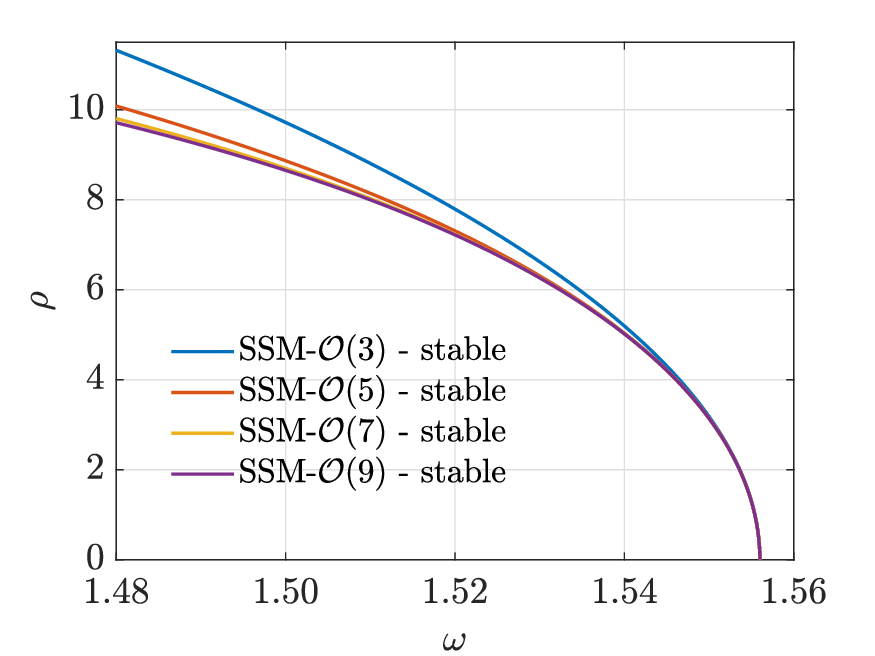} 
\includegraphics[width=.45\textwidth]{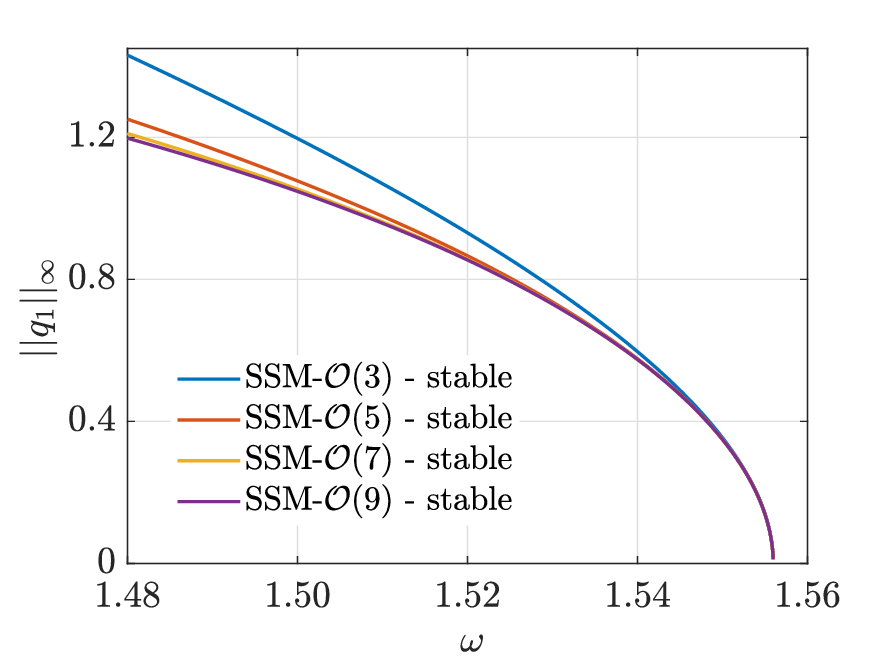}
\caption{\small Backbone curves in reduced (left panel) and physical (right panel) coordinates of the Hutchinson equation~\eqref{eq:model3_PDE}. These backbone curves are obtained via reduction on the slowest two-dimensional SSM (cf.~\eqref{eq:ROM3_stable}) at different expansion orders.}
\label{fig:mode4_backbone_unstable}
\end{figure}

When $a = \frac{\pi}{2} - 0.05 < a^* $, the first 6 eigenvalues with largest eigenvalues are described as follows
\begin{equation}
\label{eq:lamd-model3-before-hopf}
\lambda_{1,2}=-0.0230\pm1.5560\mathrm{i},\quad
\lambda_{3,4}=-0.7559\pm2.3360\mathrm{i},\quad
\lambda_{5,6}=-1.5300\pm2.7409\mathrm{i}.
\end{equation}
Since the first mode has the slowest decay rate, we take the spectral subspace of the first mode as the master subspace. The ROM on the associated SSM presents reduced dynamics below (cf.~\eqref{eq:red-auto}):
\begin{align}
\dot{\rho}& =1.906\times10^{-12}\rho^9 -2.881\times10^{-10}\rho^7 -3.888\times10^{-7}\rho^5 -0.0004353\rho^3 -0.02297\rho,\nonumber\\
\dot{\theta}& =-2.333\times10^{-11}\rho^8 -5.601\times10^{-9}\rho^6 -1.52\times10^{-6}\rho^4 -0.000593\rho^2 +1.556.
\label{eq:ROM3_stable}
\end{align}
Similar to previous examples, we obtain the damping backbone curve of the system using the ROM above. The backbone curve in both reduced and physical coordinates are shown in Fig.~\ref{fig:mode4_backbone_unstable}. Here, the backbone curve displays a softening behavior.

\begin{figure}[!ht]
\centering
\includegraphics[width=.45\textwidth]{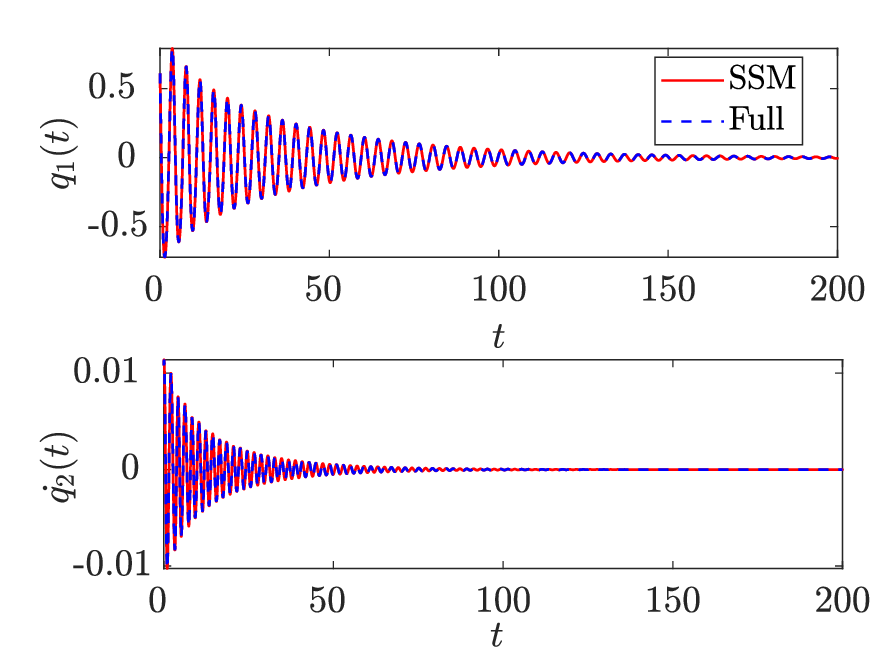}
\includegraphics[width=.45\textwidth]{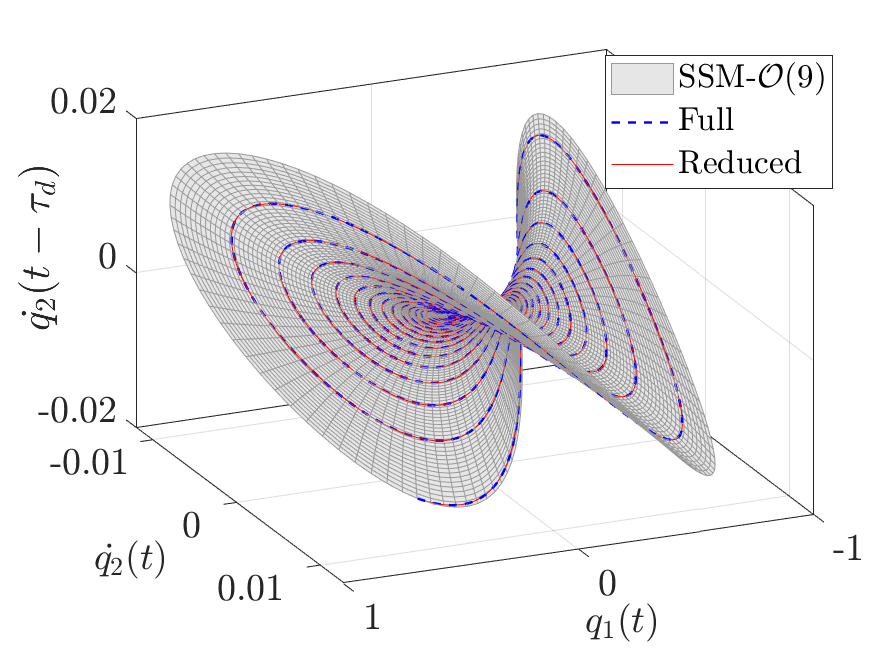} 
\caption{\small (Left panel) Time history plots of $q_1(t)$ and $\dot{q_2}(t)$ for system~\eqref{eq:model3_DDE_N}. (Right panel), the visualization of the SSM projected onto the coordinates ($q_1(t),\dot{q_2}(t),\dot{q_2}(t-\tau_d)$ (gray plane) along with the trajectories of the full system (blue dashed line) and the ROM~\eqref{eq:ROM3_stable} (red solid line).}
\label{fig:mode4_timehistory_stable}
\end{figure}

Similarly, for the $\mathcal{O}(9)$ SSM~\eqref{eq:ROM3_stable}, we use the initial condition $\bs{p}_0=( 8.1 e^{0.1\mathrm{i}}, 8.1 e^{-0.1\mathrm{i}})$. Forward-time integration is performed for both the ROM and the full system. In Fig.~\ref{fig:mode4_timehistory_stable}, the red solid line represents the SSM-based ROM simulation, while the blue dashed line represents the full system simulation. The left panel of Fig.~\ref{fig:mode4_timehistory_stable} shows that the full system's trajectory closely matches the ROM's trajectory. Additionally, the right panel projects the SSM onto the coordinates ($q_1(t),\dot{q_2}(t),\dot{q_2}(t-\tau_d)$). The numerical results demonstrate that the SSM-based ROM effectively reproduces the full system's nonlinear dynamical behavior.

\begin{figure}[!ht]
\centering
\includegraphics[width=.45\textwidth]{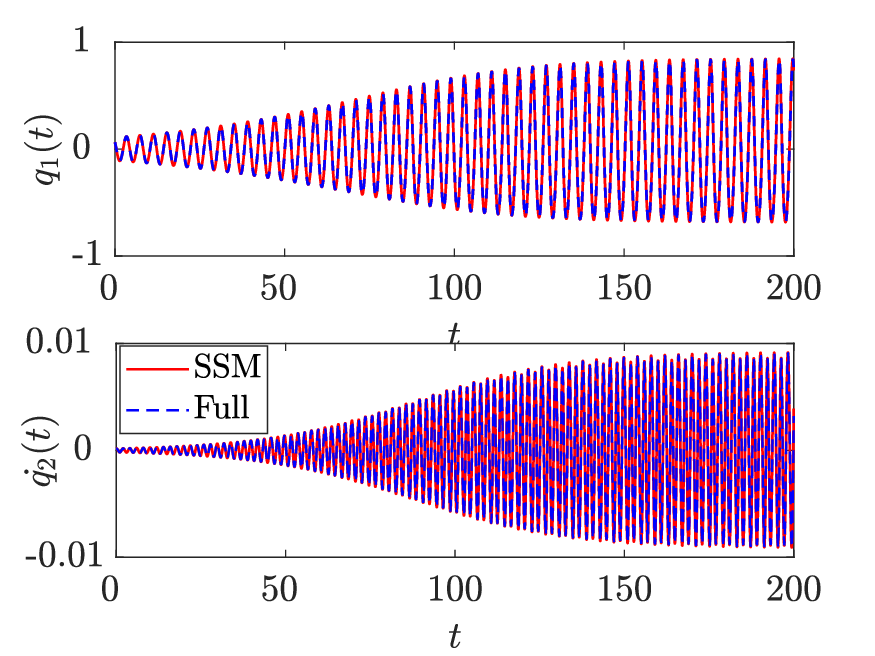}
\includegraphics[width=.45\textwidth]{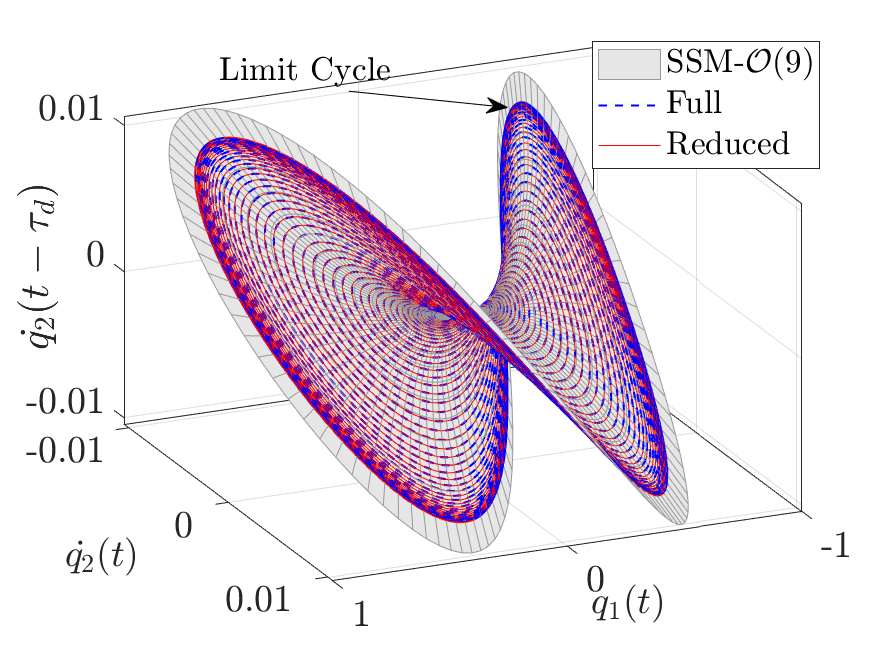} 
\caption{\small (Left panel) Time history plots of $q_1(t)$ and $\dot{q_2}(t)$ for system~\eqref{eq:model3_DDE_N}. (Right panel) The projections of the trajectories of the full system (blue dashed line) and the ROM~\eqref{eq:ROM3_unstable} (red solid line) onto the coordinates ($q_1(t), \dot{q_2}(t), \dot{q_2}(t-\tau_d)$).}
\label{fig:mode4_timehistory_unstable}
\end{figure}

\subsubsection{Reduction for dynamics after bifurcation}
Next, we focus on the dynamic behavior of system~\eqref{eq:model3_DDE_N} following a bifurcation. For this purpose, we set $a = \frac{\pi}{2} + 0.05 > a^*$. In this case, the first six eigenvalues in ascending order of their real parts are updated as 
\begin{equation}
\label{eq:lamd-model3-after-hopf}
\lambda_{1,2} = 0.0223 \pm 1.5848\mathrm{i},\quad
\lambda_{3,4} = -0.7051 \pm 2.3452\mathrm{i},\quad
\lambda_{5,6} = -1.4738 \pm 2.7436\mathrm{i}.
\end{equation}

We take the unstable subspace as the master subspace for model reduction and compute the associated two-dimensional SSM. The SSM-based ROM at $\mathcal{O}(9)$ is obtained below (cf.~\eqref{eq:red-auto}).
\begin{align}
\dot{\rho} &= a(\rho) = -4.517\times10^{-13}\rho^9 - 5.792\times10^{-10}\rho^7 - 4.028\times10^{-7}\rho^5 - 0.0004296\rho^3 + 0.02234\rho, \nonumber\\
\dot{\theta} &= b(\rho) = -1.405\times10^{-11}\rho^8 - 3.744\times10^{-9}\rho^6 - 1.146\times10^{-6}\rho^4 - 0.0005204\rho^2 + 1.585.
\label{eq:ROM3_unstable}
\end{align}

In Eq.~\eqref{eq:ROM3_unstable}, the linear term of the function $a(\rho)$ has a positive coefficient, indicating that the origin is an unstable fixed point. Further analysis of the polynomial function $a(\rho)$ reveals the existence of a non-trivial root $\rho^* = 7.0379$, satisfying $a(\rho^*) = 0$. This particular root $\rho^*$ corresponds to a limit cycle with frequency $b(\rho^*)$, which evolves from the fixed point and stays on the SSM. Numerical validation was conducted by forward time integration of the ROM and the full system under the same initial conditions: $\bs{p}_0=(1.0e^{0.1\mathrm{i}},1.0 e^{-0.1\mathrm{i}})$. In the right panel of Fig.~\ref{fig:mode4_timehistory_unstable}, we present the simulated trajectories of the full system and ROM in the coordinates ($q_1(t),\dot{q_2}(t),\dot{q_2}(t-\tau_d)$). The red solid line represents the ROM trajectory, while the blue dashed line represents the full system trajectory. These numerical results demonstrate that the SSM-based ROM effectively predicts the stable limit cycle generated by perturbations from an unstable focus and accurately predicts the response of the full system near the bifurcation point.

\section{Conclusions}
\label{sec:conclusion}
In this study, leveraging ODE approximation and spectral submanifold (SSM) theory, we have successfully developed reduced-order models (ROMs) for infinite-dimensional time-delay systems, laying a new methodological foundation for efficient and accurate analysis of high-dimensional nonlinear delay systems.

We have validated the effectiveness of this approach across three different types of systems: a delayed Duffing oscillator, two coupled oscillators with delay, and the Hutchinson equation. In each example, we first approximate the delayed system as delay-free, yet high-dimensional ODEs and further check the accuracy of the ODE approximations. We then use SSMTool~\cite{jain2021ssmtool} to construct SSM-based ROM for the approximated ODE system. We have explored the dynamical behaviors of these delay systems in the pre- and post-bifurcation regions, demonstrating that SSM-based ROMs effectively predict the nonlinear dynamics of systems, including free and forced oscillations. Particularly, SSM-based ROMs can identify isolated branches in forced response curves and predict quasi-periodic orbit bifurcations under harmonic excitation.

We have demonstrated the significant advantages of SSM-based ROMs in terms of computational efficiency and predictive accuracy via the series of progressively complex examples. Therefore, this study provides a computational framework for the efficient analysis of large-scale nonlinear delay systems in industrial applications.

\section*{Acknowledgments}
The financial support of the National Natural Science Foundation of China (No. 12302014) and Shenzhen Science and Technology Innovation Commission (No. 20231115172355001) for this work is gratefully acknowledged.

\section*{Data availability}
Data will be made available on request.

\appendix

\section{Approximating DDEs as ODEs}
In this appendix, we present a detailed analysis of the selection of discretization parameter $N$ used when the DDEs~\eqref{eq: DDEs_form} are transformed into the system of ODEs~\eqref{eq:dde-as-odes}. The selection is based on the convergence of the leading eigenvalue of the linear part of the time-delay system, and also the convergence of the nonlinear response of the DDEs~\eqref{eq:dde-as-odes} at $\epsilon=0$.

\subsection{Example~\ref{sec:ex-duffing}}
\label{sec:app-ex1}
\begin{figure}[!ht]
\centering
\includegraphics[width=.45\textwidth]{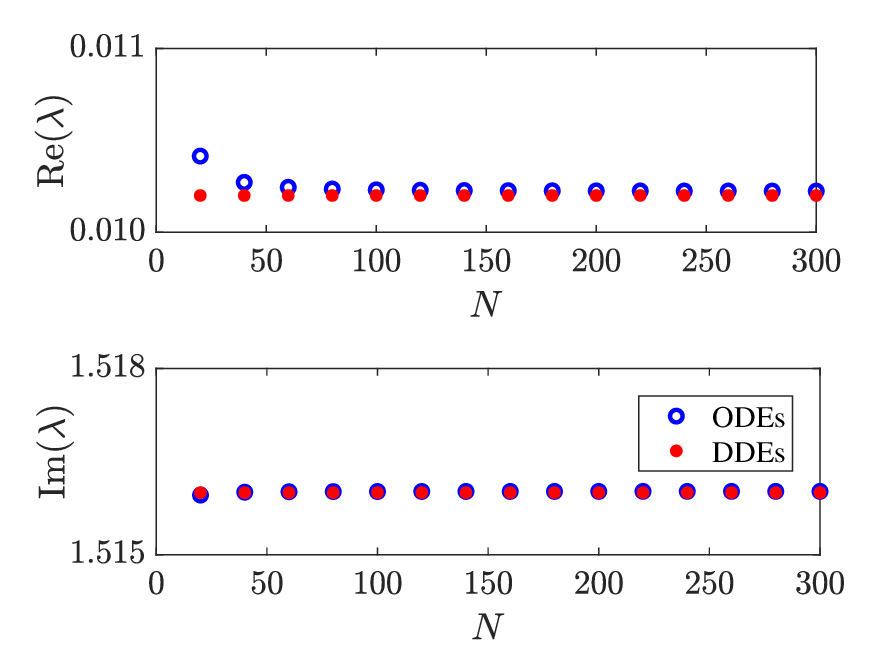} 
\includegraphics[width=.45\textwidth]{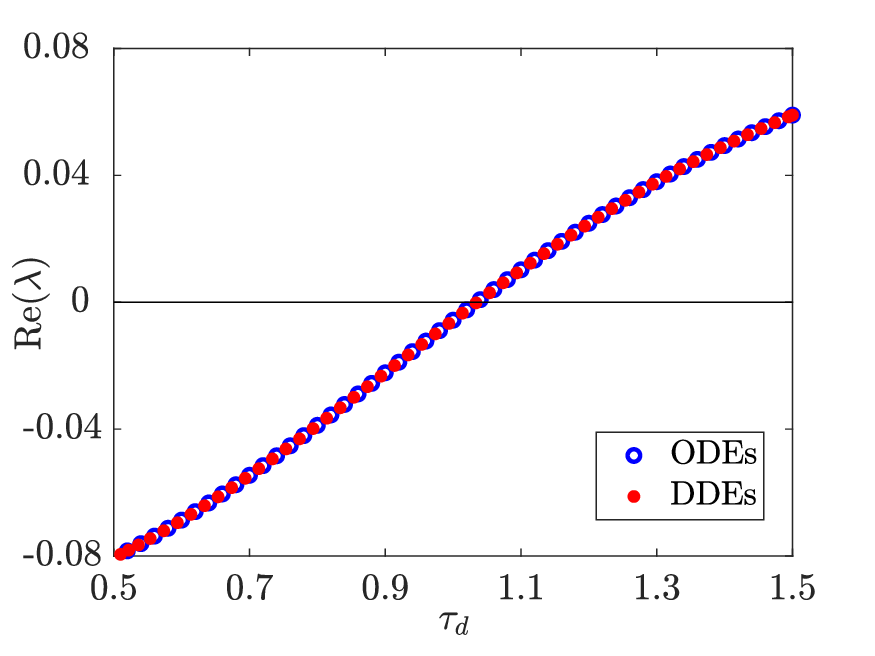}
\caption{\small The eigenvalue with the largest real part of the DDEs~\eqref{eq:duffing} and its ODE approximation system as a function of the discretization number $N$ with $\tau_d=1.1$ (left panel) and the delay parameter $\tau_d$ with $N=100$ (right panel).}
\label{fig:model1_r_distributed}
\end{figure}

We present the eigenvalue with the largest real part of the linear part of the DDEs along with its ODE approximations with various $N$ in the left panel of Fig.~\ref{fig:model1_r_distributed}. The convergence of the eigenvalue can be observed more closely from the left-upper panel of Fig.~\ref{fig:model1_r_distributed}, and the numerical results suggest that when $N$ is set to 100, the eigenvalue converges well. Furthermore, as illustrated in the right panel of Fig.~\ref{fig:model1_r_distributed}, the eigenvalue's real part increases as the delay $\tau_d$ increases. In particular, a pair of purely imaginary characteristic roots emerges upon reaching a critical value, $\tau_d=1.035$, which indicates a Hopf bifurcation in the Duffing system at $\tau_d=\tau^\ast=1.035$, the critical bifurcation point.

\begin{figure}[!ht]
\centering
\includegraphics[width=.45\textwidth]{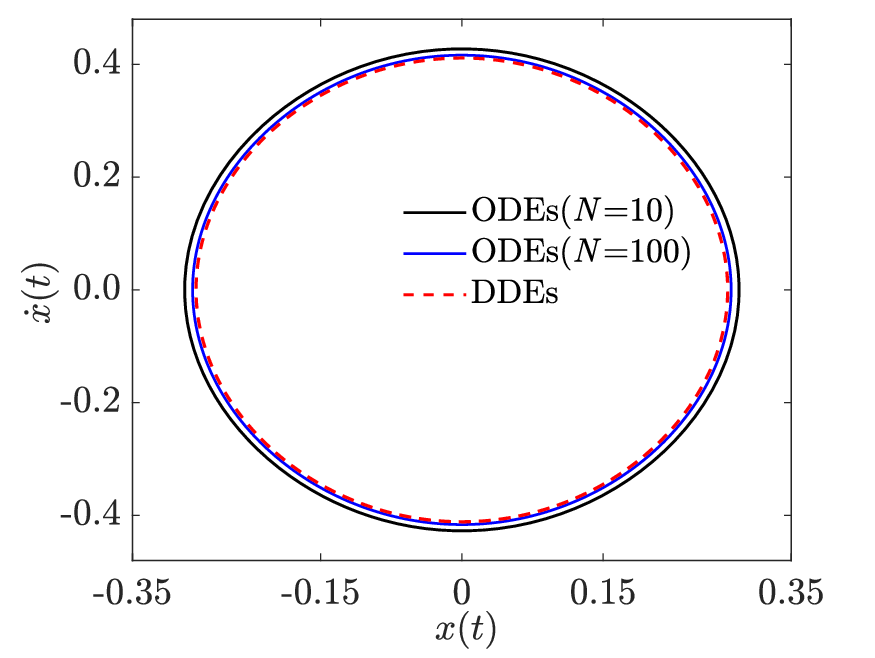} 
\caption{\small Simulated trajectories of the solutions for systems~\eqref{eq:duffing} and its ODE approximation system with $\tau_d = 1.1$.}
\label{fig:model1_solutions}
\end{figure}

We further carry out numerical analyses to compare the dynamic behaviors of the Duffing system~\eqref{eq:dde-duffing} and its corresponding ODE approximation. We let $\tau_d=1.1>\tau^\ast$ such that the system admits a limit cycle in a steady state. Fig.~\ref{fig:model1_solutions} presents the limit cycle of the original DDEs and their ODE approximations at $N=10$ and $N=100$. The solution of the original DDEs matches well with that of $N=100$, which again validates the accuracy of the ODE approximation at $N=100$. Thus, we fix $N=100$ for the ODE approximation in the rest of this example, yielding a 420-dimensional system. We then perform model reduction for this high-dimensional system using two-dimensional SSMs. 

\subsection{Example~\ref{sec:ex-two-os}}
\label{sec:app-ex2}

\begin{figure}[!ht]
\centering
\includegraphics[width=.45\textwidth]{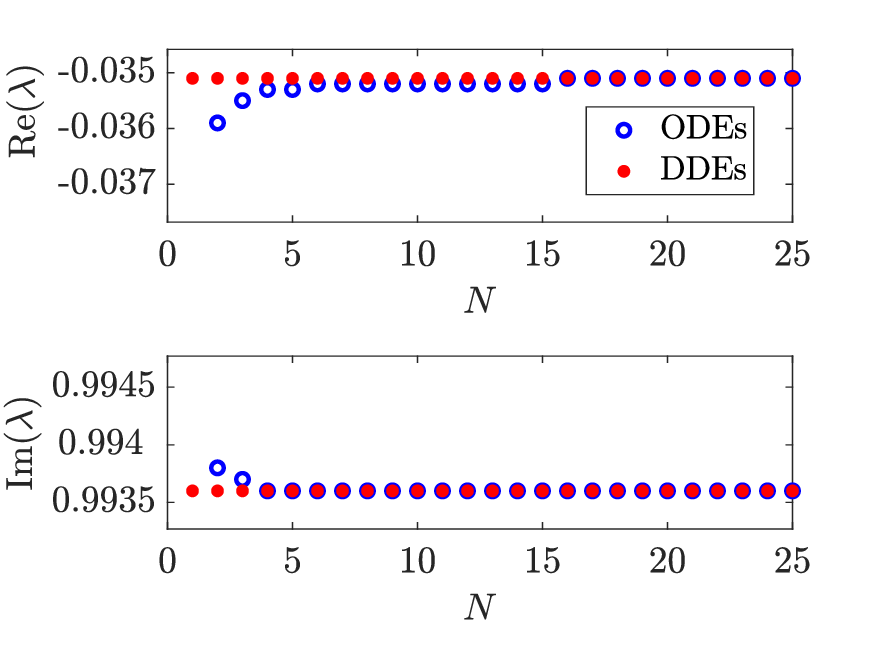} 
\includegraphics[width=.45\textwidth]{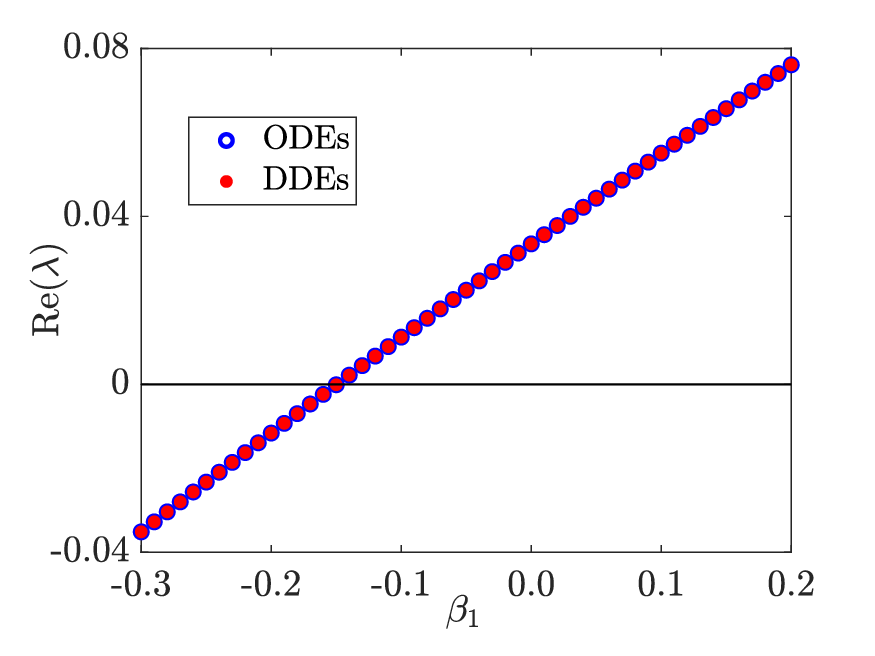}
\caption{\small The eigenvalue with the largest real part of the DDEs~\eqref{eq:model2_two_order} and its ODE approximation system as a function of the discretization number $N$ with $\beta_1=-0.3$ (left panel) and the parameter $\beta_1$ with $N=20$ (right panel).}
\label{fig:model2_r_distributed}
\end{figure}

We present the eigenvalue with the largest real part of the linear part of the DDEs along with its ODE approximation with various $N$ in the left panel of Fig.~\ref{fig:model2_r_distributed}. The convergence of the eigenvalue can be observed from the left panel of Fig.~\ref{fig:model2_r_distributed}, and the numerical results suggest that when $N$ is set to 20, the eigenvalue converges well. Furthermore, as illustrated in the right panel of Fig.~\ref{fig:model2_r_distributed}, the real part of the leading eigenvalue is increased as the parameter $\beta_1$ increases. In particular, a pair of purely imaginary characteristic roots emerges upon reaching a critical value, $\beta_1^\ast=-0.146$, which indicates a Hopf bifurcation in the system at $\beta_1=\beta_1^\ast=-0.146$.

\begin{figure}[!ht]
\centering
\includegraphics[width=.45\textwidth]{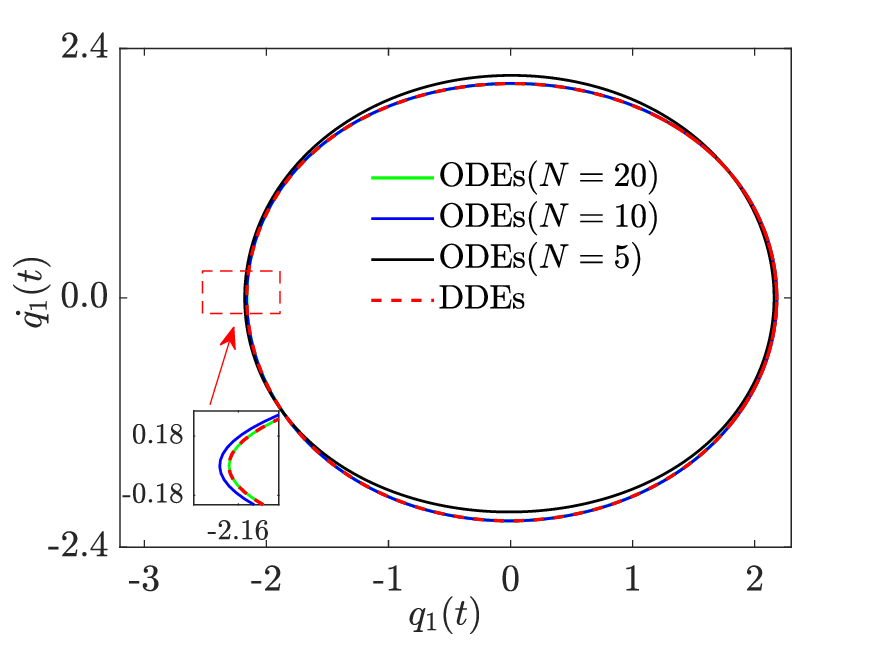} 
\caption{\small Simulated trajectories of the solutions for systems~\eqref{eq:model2_two_order} and its ODE approximation system with $\beta_1=-0.145>\beta_1^*,N=20$ and $\epsilon=0$.}
\label{fig:model2_simulation}
\end{figure}

We proceed with the simulation of the trajectories for the DDEs as denoted by Eq.~\eqref{eq:model2_two_order} and its corresponding ODE approximation system. By setting $\beta_1 = -0.145>\beta_1^* = -0.146$, the system approaches a limit cycle in the steady state. The phase portrait for the limit cycle is illustrated in Fig.~\ref{fig:model2_simulation}. In this figure, the limit cycles of the original DDEs and their ODE approximations are shown for $N = 5$, $N = 10$, and $N = 20$. The solutions of the original DDEs match well with that of the ODE approximation at $N = 20$, validating the accuracy of this approximation. Consequently, we fix $N = 20$ for the ODE approximation in the rest of this example, resulting in a 164-dimensional system. We then perform model reduction for this high-dimensional system using two-dimensional SSMs.

\subsection{Example~\ref{sec:ex-hutchinson}}
\label{sec:app-ex3}

\begin{figure}[!ht]
\centering
\includegraphics[width=.45\textwidth]{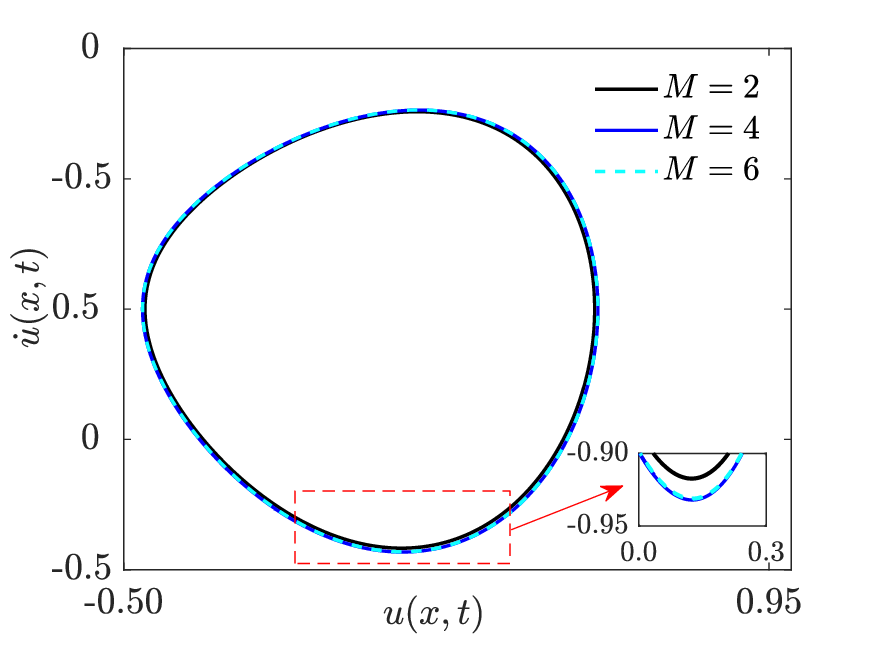} 
\caption{\small The time integration results for different $M$ in approximating the Hutchinson equation~\eqref{eq:model3_PDE} with the DDEs. Here, we take $x=1$ for visualization of the limit cycle oscillation.}
\label{fig:model3_convergence}
\end{figure}

We first determine the truncation number $M$ for which the DDE system~\eqref{eq:model3_DDE_N} can effectively approximate the Hutchinson equation~\eqref{eq:model3_PDE}. By fixing the parameter value $a = \pi/2 + 0.05$ and keeping other parameters constant, we observed that system~\eqref{eq:model3_PDE} admits a limit cycle oscillation because the system undergoes a supercritical Hopf bifurcation. We choose $M = 2, 4, 6$ and apply forward simulation via dde23 of \textsc{matlab} to the DDE system~\eqref{eq:model3_DDE_N}. The produced limit cycles in steady state under various $M$ are shown in Fig.~\ref{fig:model3_convergence}, from which we see that the DDEs~\eqref{eq:model3_DDE_N} with $M=4$ can effectively approximate the Hutchinson equation~\eqref{eq:model3_PDE}. We select $M=4$ for the remaining computations.

\begin{figure}[!ht]
\centering
\includegraphics[width=.45\textwidth]{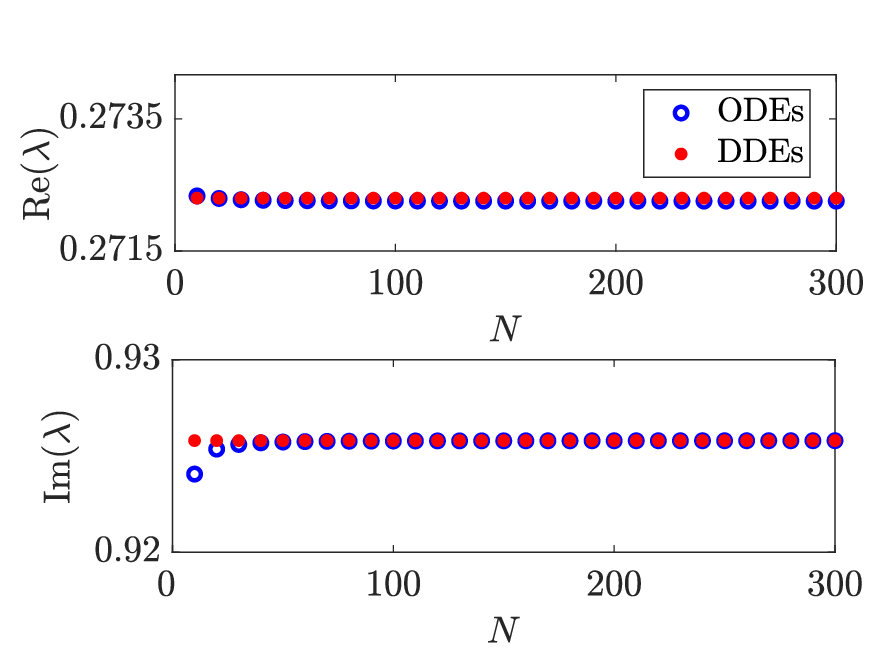} 
\includegraphics[width=.45\textwidth]{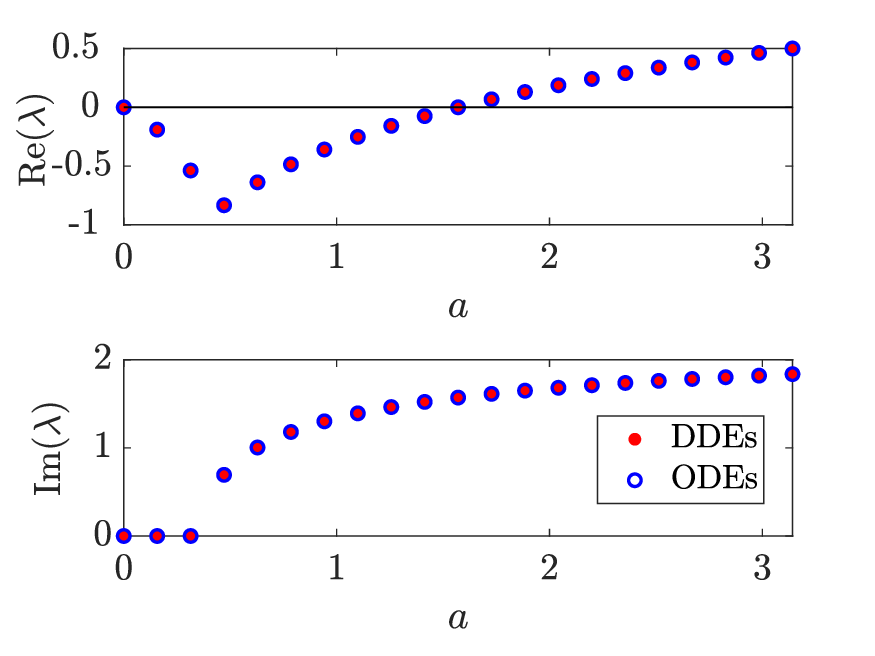} 
\caption{\small The eigenvalue with the largest real part of the DDEs~\eqref{eq:model3_DDE_N} and its ODE approximation system as a function of the discretization number $N$ with $a=\frac{\pi}{2}+0.05$ (left panel) and the delay parameter $a$ with $N=100$ (right panel).}
\label{fig:model3_r_distributed}
\end{figure}

Next, we again follow Section~\ref{sec:dde-to-ode} to transform the DDEs~\eqref{eq:model3_DDE_N} into a system of ODEs~\eqref{eq:non_eq1}. Here, we have $n=4$, and hence $\bs{z}\in\mathbb{R}^{4(2N+1)}$. Similar to the previous examples, we present the leading eigenvalue of the linear part of the DDEs, along with its ODE approximations with various $N$ in the left panel of Fig.~\ref{fig:model3_r_distributed}. From the left panel, we observe that the eigenvalue converges well at $N=100$. Furthermore, as illustrated in the right panel of Fig.~\ref{fig:model3_r_distributed}, the real part of the leading eigenvalue passes through $a=\pi/2$, consistent with the analytic prediction given by~\eqref{eq:model3_lambda_eq}. 

\begin{figure}[!ht]
\centering
\includegraphics[width=.45\textwidth]{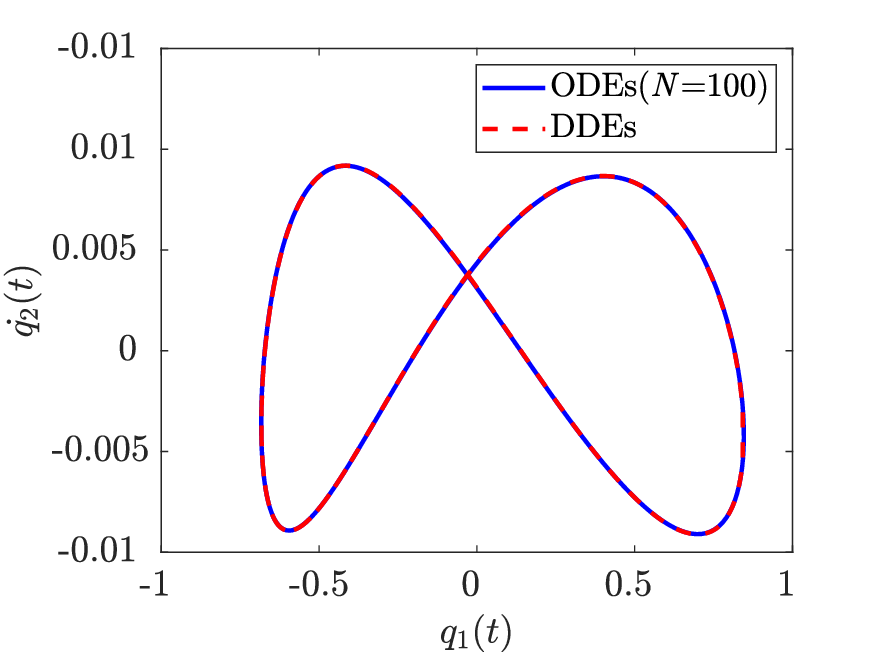} 
\caption{\small Simulated limit cycle for the DDE system (\ref{eq:model3_DDE_N}) with $a=\frac{\pi}{2}+0.05$ and its ODE approximation system.}
\label{fig:model3_simulation}
\end{figure}

To further confirm that $N=100$ is sufficient for the ODE approximation, a forward simulation was conducted. The simulated results are visualized in Fig.~\ref{fig:model3_simulation}, with the blue dashed line representing the trajectory of the ODE approximation system and the red curve representing the trajectory of the DDEs. The results show a good match between the trajectories of the two systems, validating the approximation system's accuracy when $N$ = 100.

\section{Transformation of initial conditions}
\label{sec:transform-ini}

Here we present how to determine the initial conditions of the SSM-based ROM~\eqref{eq:red-auto} from the initial conditions of the DDEs~\eqref{eq: DDEs_form}. We transform the initial time history~\eqref{eq: DDEs_initi} to the initial conditions of the SSM-based ROM~\eqref{eq:red-auto} via two steps below:
\begin{itemize}
    \item \emph{Transforming initial conditions of DDEs~\eqref{eq: DDEs_form} to the associated ODEs~\eqref{eq:non_eq1}}. We recall that $\bs{z}$ in~\eqref{eq:non_eq1} is defined as \( \bs{z} = (\bs{u}_0, \dots, \bs{u}_N, \bs{w}_1, \dots, \bs{w}_N) \). Let $\tau_i=i\tau_\mathrm{d}/N$, we follow from equations~\eqref{eq:DDEs_transform}--\eqref{eq:dde-as-odes} that $\bs{u}_i(0) = \bs{x}(-\tau_i)$ for $0\leq i\leq N$ and $\bs{w}_i(0) = \dot{\bs{x}}(-\tau_i)$ for $1\leq i\leq N$.
   We further know from~\eqref{eq: DDEs_initi} that $\bs{x}(-\tau_i)=\bs{x}_0(-\tau_i)$ and $\dot{\bs{x}}(-\tau_i)=\dot{\bs{x}}_0(-\tau_i)$. Thus, we have obtained $\bs{z}(0)$ from the initial time history $\bs{x}_0$ shown in~\eqref{eq: DDEs_initi}.
    \item \emph{Projection of initial conditions to reduced coordinates}. Recall that $\bs{p}=(p,\bar{p})$ is the vector of reduced coordinates. We simply take a linear projection $p(0)=\left(\bs{v}^\mathcal{E}\right)^\mathrm{T}\bs{z}_0$ to determine the initial reduced coordinates, where $\bs{v}^\mathcal{E}$ is the eigenvector associated with the master subspace $\mathcal{E}$. Then we determine the initial conditions for the SSM-based ROM~\eqref{eq:red-auto} via $p=\rho e^{\mathrm{i}\theta}$. Here we take the linear projection because it is simple to implement and provides good approximation in practice. To improve accuracy, one may determine $\bs{p}(0)$ by minimizing the distance from $\bs{z}_0$ to the SSM, namely, $\bs{p}(0)=\mathrm{argmin}_{\bs{p}}||\bs{z}_0-\bs{W}(\bs{p})||$.
\end{itemize}

\section{Supplementary analysis to example~\ref{sec:ex-duffing}} \label{appendix:C}
We revisit Section~\ref{subsec:5.1.2} but take $\tau_d=1.75$ instead of $\tau_d=1.1$. We follow the same procedures as in Section~\ref{subsec:5.1.2} to perform SSM-based reduction. Here, the first five eigenvalues of the linear portion of the approximated system, corresponding to~\eqref{eq:dde-duffing}, are:

\begin{align}
\lambda_{1,2} =  0.0756 \pm 1.4614\mathrm{i}, \quad
\lambda_{3} = -1.5153 \pm  0.0000\mathrm{i}, \quad
\lambda_{4,5} = -1.7444 \pm 4.2204\mathrm{i}.
\label{eq:lamd-model1-post-hopf}
\end{align}

We compute the SSM corresponding to the first pair of eigenvalues and formulate an ROM on the SSM. In the scenario of free vibration where $\epsilon=0$, we derive a two-dimensional ROM~\eqref{eq:ROM1_post_hopf_fail} for the SSM approximated at $\mathcal{O}(15)$. The ROM is listed as follows:

\begin{align}
\label{eq:ROM1_post_hopf_fail}
\dot{\rho} &= -6.711 \times 10^{-16} \rho^{15} - 5.134 \times 10^{-14} \rho^{13} - 4.092 \times 10^{-12} \rho^{11} - 3.437 \times 10^{-10} \rho^9 \nonumber\\
&\quad - 3.104 \times 10^{-8} \rho^7 - 3.243 \times 10^{-6} \rho^5 - 0.0002964 \rho^3 + 0.07559 \rho, \nonumber\\
\dot{\theta} &= -4.193 \times 10^{-15} \rho^{14} - 3.205 \times 10^{-13} \rho^{12} - 2.555 \times 10^{-11} \rho^{10} - 2.165 \times 10^{-9} \rho^8 \nonumber\\
&\quad - 2.024 \times 10^{-7} \rho^6 - 2.273 \times 10^{-5} \rho^4 - 0.004277 \rho^2 + 1.461.
\end{align}

\begin{figure}[!ht]
 \centering 
 \includegraphics[width=0.45\textwidth]{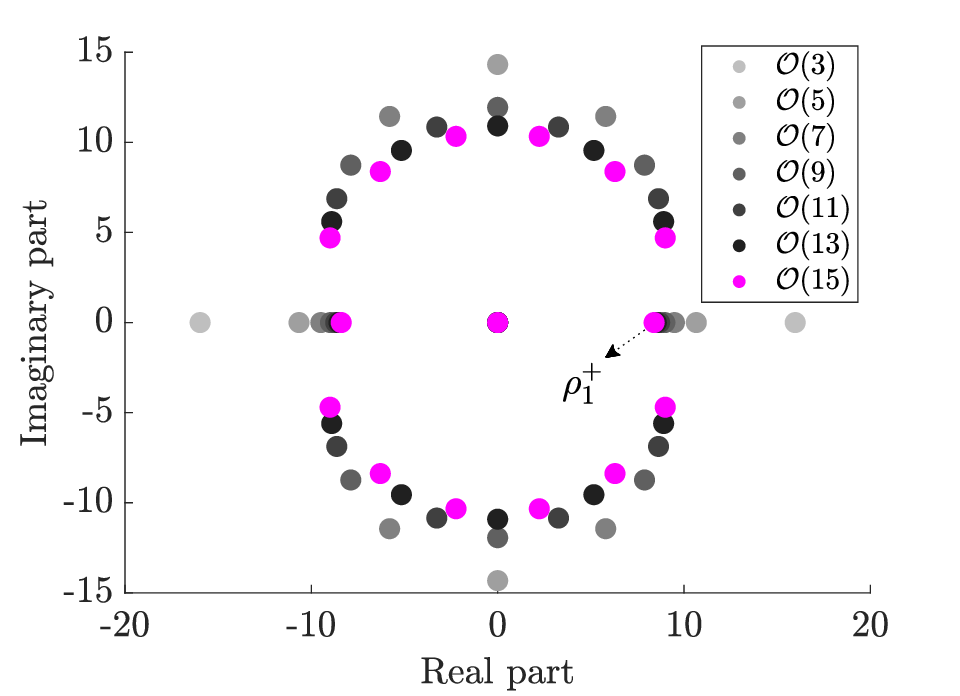} 
  \includegraphics[width=0.45\textwidth]{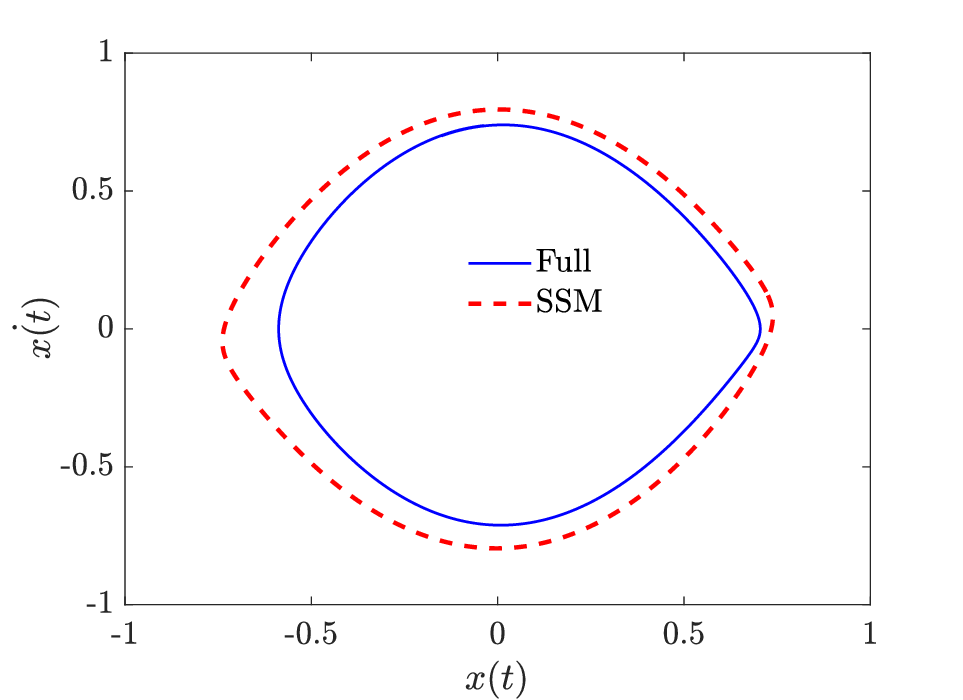} 
 \caption{\small (Left panel) Roots of \( a(\rho) \) in~\eqref{eq:ROM1_post_hopf_fail} computed using SSMTool with truncation orders up to \( \mathcal{O}(15) \). Darker markers indicate higher-order approximations, while magenta markers highlight zeros from the highest-order truncation. (Right panel) Predicted limit cycle via the spurious root \( \rho_1^+ \) along with the limit cycle obtained from the full system (cf. Eq.~\eqref{eq:dde-as-odes}).}
 \label{fig:root_plot_fail}
\end{figure}

At $\tau_d=1.75$, the system also admits a limit cycle in this post-flutter region, as seen in the right panel of Fig.\ref{fig:root_plot_fail}. This limit cycle has a larger radius because $\tau_d=1.75$ moves further away from $\tau^\ast=1.035$. The SSM-based ROM~\eqref{eq:ROM1_post_hopf_fail} fails to predict this limit cycle. In particular, we compute the roots of $a(\rho)$ in~\eqref{eq:ROM1_post_hopf_fail} with increasing truncation orders and obtain the results shown in the left panel of Fig.\ref{fig:root_plot_fail}. We observe that $a(\rho)$ does not admit non-spurious, non-trivial roots, which is different from the case of Fig.\ref{fig:model1_roots}. Thus, the SSM-based ROM fails to predict the existence of limit cycles. We note that the real root $\rho_1^+$ in Fig.\ref{fig:root_plot_fail} is on the boundary of the domain of convergence, and hence this root is spurious f~\cite{ponsioen2019analytic}. Indeed, the limit cycle corresponding $\rho_1^+$ deviates significantly from the reference solution, as seen in the right panel of Fig.\ref{fig:root_plot_fail}.

\begin{figure}[!ht]
 \centering 
  \includegraphics[width=0.45\textwidth]{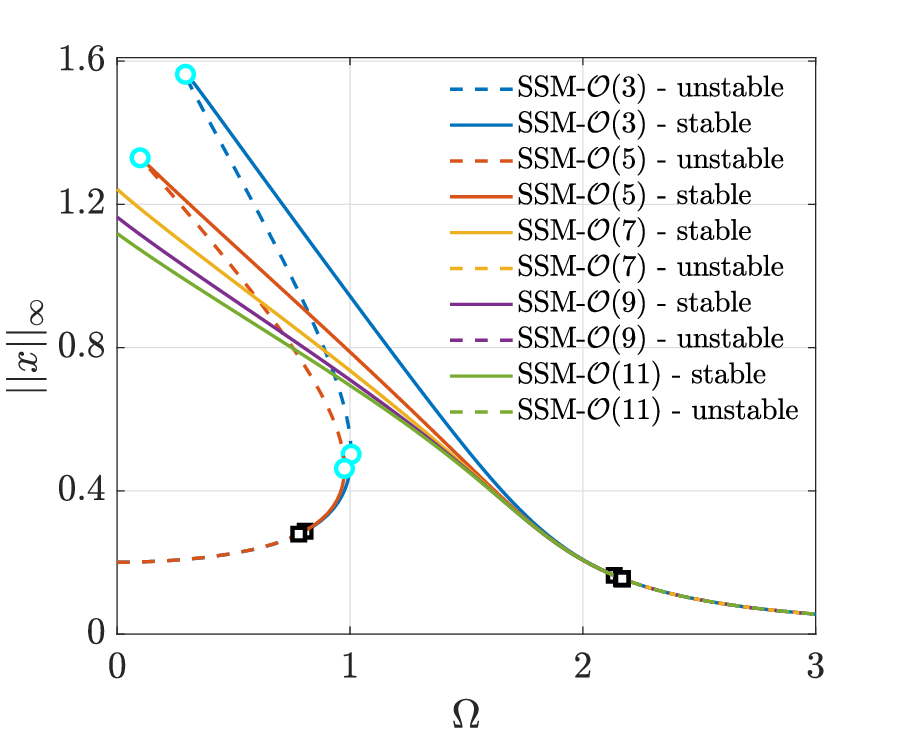} 
 \caption{\small Forced response curve (FRC) for periodic orbits of the delayed Duffing system (Eq.~\eqref{eq:dde-duffing}) at excitation amplitude \( \epsilon = 0.4 \). Solid and dashed lines represent stable and unstable periodic orbits, respectively. Cyan circles indicate saddle-node bifurcations, while black squares mark Hopf bifurcations (HB).}
 \label{fig:frc-plot-fail}
\end{figure}

We further revisit Section \ref{sec:5.1.3} but increase $\epsilon=0.4$. The obtained FRCs at various expansion orders are presented in Fig.~\ref{fig:frc-plot-fail}. We observe that the predicted forced periodic orbits fail to converge if the response has a large amplitude. This is expected because we have adopted a leading-order approximation to the non-autonomous part of SSMs. Such an approximation is not sufficient when $\epsilon$ is large, as discussed in Remark~\ref{rmk:ssm-reduction}.

\end{document}